\documentclass{amsart}%
\usepackage{eurosym}
\usepackage{amsmath}
\usepackage{amsfonts}
\usepackage{amssymb}
\usepackage{graphicx}
\usepackage{latexsym}
\usepackage{amscd}
\usepackage{tikz-cd}
\usepackage{enumitem}
\usepackage[all, cmtip]{xy}
\usepackage[utf8]{inputenc}
\usepackage{amsthm}
\usepackage{xcolor}
\usepackage[a4paper,  margin=1in]{geometry}
\usepackage{epic}
\usepackage{pstricks}
\usepackage[all]{xy}
\usepackage{graphicx}
\usepackage{epic,eepic}
\usepackage{epsfig}
\usepackage{float}
\usepackage{extarrows}
\usepackage{tikz}
\usepackage{calc}%
\usepackage{mathtools}
\providecommand{\U}[1]{\protect\rule{.1in}{.1in}}
\providecommand{\U}[1]{\protect\rule{.1in}{.1in}}
\providecommand{\U}[1]{\protect\rule{.1in}{.1in}}
\providecommand{\U}[1]{\protect\rule{.1in}{.1in}}
\newtheorem{theorem}{Theorem}[section]

\newtheorem{conjecture}[theorem]{Conjecture}
\newtheorem{corollary}[theorem]{Corollary}

\newtheorem{definition}[theorem]{Definition}
\newtheorem{example}[theorem]{Example}

\newtheorem{lemma}[theorem]{Lemma}

\newtheorem{proposition}[theorem]{Proposition}
\newtheorem{remark}[theorem]{Remark}

\newtheorem{convention}[theorem]{Convention}

\def\id{\mathop{\text{Id}}}

\def\remove#1{}

\def\path{\mathop{\hbox{\rm Path}}}

\newcommand{\K}{{\mathbb{K}}}

\DeclareMathOperator{\Span}{span}
\DeclareMathOperator{\End}{End}

\newcommand{\Rc}{\check{R}}
\DeclareMathOperator{\tr}{tr}
\newcommand{\kk}{\Bbbk}
\DeclareMathOperator{\vect}{vec}
\DeclareMathOperator{\rank}{rank}
\newcommand{\Qhat}{\widehat{Q}}

\usetikzlibrary{decorations.markings}

\tikzstyle{vertex}=[circle, draw, fill, inner sep=0pt, minimum size=6pt]
\begin{document}
\title[Quivers with
quantum Yang--Baxter equation and Hecke condition]{Quivers with
quantum Yang--Baxter equation and Hecke condition: Deformation of face algebras}
\author{Cody Gilbert}
\address{Department of Mathematics and Statistics, Saint Louis University, St. Louis,
MO-63103, USA}
\email{cody.gilbert@slu.edu}
\author{Ashish K. Srivastava}
\thanks{Ashish K Srivastava is the corresponding author}
\address{Department of Mathematics and Statistics, Saint Louis University, St. Louis,
MO-63103, USA}
\email{ashish.srivastava@slu.edu}
\dedicatory{Dedicated to the memory of Professor C. N. Yang}

\vspace{150mm}

\begin{abstract}
In this paper, we initiate the study of quivers carrying quantum Yang--Baxter and Hecke
structure. Calling a quiver $Q$ to be a solution of the QYBE when the adjacency
matrix of its Kronecker square is a solution, we show that this holds precisely
when the adjacency matrix $A$ satisfies $A^2 = \mu A$ for a scalar $\mu$. We determine exactly when
Kulish's rank-one construction yields a Hecke $R$-matrix that satisfies both the braided QYBE and the
Hecke condition. We deform Hayashi's face algebra by the resulting RTT
relations, and prove that the quantum matrix algebra $\mathcal{O}_q(M_n)$
is isomorphic as a bialgebra to the Hecke-deformed face algebra of the rose
quiver with $n$ petals, and that for a disjoint union of $m$ such roses the
deformation is a genuine quantum groupoid with $m$-dimensional base,
isomorphic as an algebra to $m^2$ copies of $\mathcal{O}_q(M_n)$.

\bigskip


\bigskip
	
\end{abstract}
\maketitle

\bigskip

\section{Introduction}

\bigskip

\noindent The quantum Yang-Baxter Equation (QYBE) is a fundamental consistency condition that arises in mathematical physics, quantum groups, and noncommutative algebra. Its origin traces back to influential papers of Yang \cite{Yang} on quantum mechanical many-body problem and Baxter \cite{Baxter} on solution of the eight vertex model. A solution to the quantum Yang-Baxter equation is a linear map $R: V\otimes V\rightarrow V\otimes V$ where $V$ is a vector space over a field such that the equality $(R\otimes I_V)(I_V \otimes R)(R\otimes I_V)=(I_V\otimes R)(R\otimes I_V)(I_V\otimes R)$ of operators holds in the space $\End(V\otimes V\otimes V)$. The matrix of such a linear map is called an $R$-matrix. The term quantum Yang–Baxter Equation (QYBE) was coined by Faddeev in the late 1970s. Intuitively, the QYBE ensures that when three interacting particles or algebraic objects are rearranged, the final outcome does not depend on the order in which the pairwise interactions occur. It has applications in the study of exactly solvable models in statistical mechanics and quantum field theory, and has connections to knot theory, braid theory and Hopf algebras. In \cite{D} Drinfeld suggested to study the set-theoretical solutions of QYBE. A set-theoretical solution is a solution for which $V$ is a vector space spanned by a set $X$, and $R$ is the linear operator induced by a mapping $X\times X \rightarrow X\times X$. For geometric and algebraic interpretations of such solutions, the reader is referred to Etingof et al. \cite{ESS}. 

Closely linked to the QYBE is the Hecke condition, which places an additional polynomial constraint on the $R$-matrix. The Hecke condition states that $(R-qI)(R+q^{-1}I)=0$ for some nonzero parameter $q$. Together, the QYBE and the Hecke condition provide a powerful algebraic framework used to construct quantum symmetries and invariants. Solutions that satisfy both conditions lead to braided monoidal categories, quantum knot invariants like the Jones polynomial, and integrable lattice models in statistical mechanics. 

In this paper we develop a quiver-theoretic perspective on this framework. We say that a quiver $Q$ satisfies the QYBE if the adjacency matrix of the Kronecker square quiver satisfies the QYBE. We show that if a quiver $Q$ satisfies the QYBE then its adjacency matrix $A$ satisfies the ``Temperley-Lieb idempotent" property, that is, $A^2=\mu A$ for some scalar $\mu$. In particular, we show that a finite symmetric quiver $Q$ satisfies the QYBE if and only if its adjacency matrix is a scalar multiple of a symmetric idempotent; for nonnegative real (e.g. genuine weighted) adjacency matrices this means each nontrivial component is a complete weighted quiver.

In Section 3 of this paper, we give a recipe to construct Hecke $R$-matrices from the weighted adjacency matrices of such complete quivers via Kulish's rank-one tensor construction, and we study the four-cycle identity that characterizes the resulting weighted quivers. 

In Section 4, we study Hecke deformed face algebras, where local RTT relations at vertices with non-trivial loop space deform the face algebra structure. Hayashi introduced face algebras as a new class of quantum groups inspired by the quantum inverse scattering method and solvable lattice models of face type. Face algebras also arise in Jones' index theory of subfactors. Face algebras are related to Fusion rules of Wess-Zumino-Witten model in Conformal Field Theory. The class of Face algebras contains all bialgebras and just like bialgebras, face algebras also produce monoidal categories as their (co-)module categories. Recently Huang, Walton, Wicks, and Won have shown that the weak bialgebras that coact universally on the path algebra $\Bbbk Q$ (either from the left, from the right, or from both directions compatibly) are each isomorphic to Hayashi's face algebra $\mathfrak H(Q)$. In \cite[Conjecture 1.6]{HWWW} it is conjectured that if $H$ is a finite-dimensional weak bialgebra over an algebraically closed field $\Bbbk$ with commutative counital subalgebras, then $H$ is isomorphic to a weak bialgebra quotient of face algebra $\mathfrak H(Q)$ for some finite quiver $Q$. This motivates us to study quotients of face algebra. We deform Hayashi's face algebra by the resulting RTT relations, and prove that the quantum matrix algebra $\mathcal{O}_q(M_n)$ is isomorphic as a bialgebra to the Hecke-deformed face algebra of the rose quiver with $n$ petals, and that for a disjoint union of $m$ such roses the deformation is a genuine quantum groupoid with $m$-dimensional base, isomorphic as an algebra to $m^2$ copies of $\mathcal{O}_q(M_n)$.

\section{Quivers Satisfying the QYBE}

\begin{definition}
    Let $Q$ be a quiver (that is, a directed graph) with the adjacency matrix $A=[a_{ij}]$. The Kronecker square $\widehat{Q}$ is a quiver whose adjacency matrix is $A\otimes A$ where $A\otimes A$ is the Kronecker product matrix given as $A\otimes A=[b_{ij}]$ with $b_{ij}=a_{ij}A$. 
\end{definition} 
Let $Q=(Q_0, Q_1, s, r)$ be a quiver where $Q_0$ is the set of vertices, $Q_1$ is the set of arrows and both $s$, and $r$ are functions from $Q_1$ to $Q_0$ such that for any arrow $e$ from $u$ to $v$, $s(e)=u$ and $r(e)=v$.
    Then the Kronecker square $\widehat{Q}=(\widehat{Q_0}, \widehat{Q_1}, \widehat{s}, \widehat{r})$ is described as
    \begin{align*}
        &\widehat{Q_0} = \{[v_i, v_j]\, : \, v_i, v_j\in Q_0\}\\
        &\widehat{Q_1} = \{[e_i, e_j]\, : \, e_i, e_j\in Q_1\}\\
        &\qquad\widehat{s}([e,f])=[s(e), s(f)]\\
        &\qquad \widehat{r}([e,f])=[r(e), r(f)]
    \end{align*}

Recall that an $n^2\times n^2$ matrix $X$ is said to satisfy the quantum Yang-Baxter Equation (QYBE) if $$(X\otimes I_n)(I_n\otimes X)(X\otimes I_n)=(I_n\otimes X)(X\otimes I_n)(I_n\otimes X)$$ where $I_n$ is the $n\times n$ identity matrix. 

\begin{definition}
Let $Q$ be a quiver with adjacency matrix $A$. We will say that $Q$ satisfies the quantum Yang-Baxter Equation if the adjacency matrix $X=A\otimes A$ of $\widehat{Q}$ satisfies the quantum Yang-Baxter Equation. 
\end{definition}

We first investigate when an adjacency matrix $X=A\otimes A$ satisfies the QYBE. It will quickly be apparent that such $X$ have a highly degenerate structure. This means $X$ will generally not be an $R$-matrix as it will be singular (barring the identity matrix); however, in the next section, we will use such $X$ to produce Hecke $R$-matrices, $\check{R}(q)$ which arise in a straightforward, combinatorial way via weighted adjacency matrices of complete quivers.

A quick computation introduces the equation of interest for this section. Suppose we have a matrix $Y=A\otimes B\in \mathbb{C}^{n^2\times n^2}$ where $A,B\in C^{n\times n}$. We want to determine when $Y$ satisfies the QYBE. We must have 
\begin{align*}
    &\qquad \qquad \qquad \qquad \,\,\,\,\,(Y\otimes I_n)(I_n\otimes Y)(Y\otimes I_n)=(I_n\otimes Y)(Y\otimes I_n)(I_n\otimes Y)\\
    &\iff (A\otimes B\otimes I_n)(I_n\otimes A\otimes B)(A\otimes B\otimes I_n)=(I_n\otimes A\otimes B)(A\otimes B\otimes I_n)(I_n\otimes A\otimes B)\\
    &\iff\qquad \qquad \quad\, (A\otimes BA\otimes B)(A\otimes B\otimes I_n)=(A\otimes AB\otimes B)(I_n\otimes A\otimes B)\\
    &\iff \qquad \qquad \qquad\qquad \qquad \,\,\,A^2\otimes BAB\otimes B=A\otimes ABA\otimes B^2
\end{align*}

\noindent In the case where $X=A\otimes A$, the above equation becomes 
\begin{equation}\label{eqn:Kron-QYBE}
A^2\otimes A^3\otimes A=A\otimes A^3\otimes A^2.
\end{equation} 

The main theorem of this section, Theorem \ref{thm:R-Kronecker}, characterizes the matrices $A$ that satisfy equation (\ref{eqn:Kron-QYBE}) above and finds that such matrices have the ``Temperley-Lieb idempotent" property, that is, $A^2=\mu A$ for some scalar $\mu$.

\begin{theorem}\label{thm:R-Kronecker}
    Let $A$ be a $n\times n$ matrix, such that $A^3\neq 0$, then $$A\otimes A^3\otimes A^2=A^2\otimes A^3\otimes A$$ if and only if $A^2\otimes A=A\otimes A^2$ if and only if $A^2=\mu A$.
\end{theorem}

\noindent Theorem \ref{thm:R-Kronecker} is an immediate consequence of Lemma \ref{lemma:scalar_prop}, Lemma \ref{lemma:cubic->square} and Lemma \ref{lemma:square->scalar} below.

\begin{lemma}\label{lemma:scalar_prop}
    Suppose $X, Y, U, V$ are nonzero square matrices of the same size, then if $X\otimes Y=U\otimes V$, then there exists a nonzero scalar $\lambda$ such that $X=\lambda U$ and $Y=\frac{1}{\lambda}V$.
\end{lemma}
\begin{proof}
Let $X,Y,U,V$ be $n\times n$ matrices over $\mathbb{R}$ (a similar argument works for $\mathbb{C}$). Define $$x=\text{vec}(X), \quad y=\text{vec}(Y), \quad u=\text{vec}(U), \quad v=\text{vec}(V).$$ The equality $X\otimes Y=U\otimes V$ thus corresponds to the equality $x\otimes y=u\otimes v$. Note, this is not to say $\text{vec}(X\otimes Y)=x\otimes y$, but rather, if $x\otimes y=u\otimes v$, then the entries of $X\otimes Y$ must be the same as the entries of $U\otimes V$, so it suffices to show $x\otimes y=u\otimes v$. Furthermore, the equality $x\otimes y=u\otimes v$ is equivalent to the equality $xy^T=uv^T$ of rank one matrices. 

As a rank one matrix, the singular value decomposition of $xy^T$ can be written as $xy^T=\sigma pq^T$, where $$\sigma=||x||\cdot ||y||,\quad  p=\frac{x}{||x||}, \quad q=\frac{y}{||y||}.$$ Similarly, one can write $uv^T=\tau rs^T$ where $$\tau=||u||\cdot||v||,\quad r=\frac{u}{||u||},\quad s=\frac{v}{||v||}.$$ As $xy^T=uv^T$, we must have $\sigma pq^T=\tau rs^T$. 

For a nonzero rank one matrix, the left and right singular values are unique up to sign, thus there is a scalar $\alpha$ with $$p=\alpha r,\quad q=\frac{1}{\alpha}s,\quad \sigma=\tau.$$ Substituting into the equations $p=\frac{x}{||x||}$ and $q=\frac{y}{||y||}$ we find $$x=\left(\frac{\alpha||x||u}{||u||}\right),\qquad y=\left( \frac{\alpha^{-1}||y||v}{||v||} \right).$$ 

\noindent Letting $\lambda:= \frac{\alpha||x||}{||u||}$, it follows $x=\lambda u$ and, as $\tau=||u||\cdot ||v||=||x||\cdot ||y||=\sigma$, $$\lambda^{-1}v=\frac{\alpha^{-1}||u||}{||x||}v=\frac{\alpha^{-1}||y||}{||v||}v=y.$$

\noindent As $x=\lambda u$ and $y=\lambda^{-1}v$, we must have $X=\lambda U$ and $Y=\frac{1}{\lambda}V$, as desired.
\end{proof}

\begin{lemma}\label{lemma:cubic->square}
    Let $A$ be an $n\times n$ matrix and suppose $A^3\neq 0$, then $A^2\otimes A=A\otimes A^2$ if and only if $A^2\otimes A^3\otimes A=A\otimes A^3\otimes A^2$.
\end{lemma}
\begin{proof}
    Suppose $A^2\otimes A=A\otimes A^2$. As $A^2\otimes A=A\otimes A^2$, we must have
    \begin{align*}
        (A\otimes A)(A^2\otimes A)&=(A\otimes A)(A\otimes A^2)\\
       \Rightarrow A^3\otimes A^2&=A^2\otimes A^3
    \end{align*}

\noindent Furthermore, if $I$ is the $n\times n$ identity matrix, then we find 
\begin{align*}
    (A\otimes I)(A^2\otimes A)&=(A\otimes I)(A\otimes A^2)\\
    \Rightarrow A^3\otimes A&=A^2\otimes A^2
\end{align*}
\noindent and similarly, $A^2\otimes A^2=A\otimes A^3$. All together, we have the equality $$A\otimes A^3=A^2\otimes A^2=A^3\otimes A.$$ Using the above equalities, we get the equation $$A^2\otimes A^3\otimes A=A^3\otimes A^2\otimes A=A^3\otimes A\otimes A^2=A\otimes A^3\otimes A^2.$$ 

For the converse, we write $$x=\text{vec}(X), \qquad y=\text{vec}(A^3)\qquad z=\text{vec}(A^2).$$ The equality $A\otimes A^3\otimes A^2=A^2\otimes A^3\otimes A$ is equivalent to the equality $x\otimes y\otimes z=z\otimes y\otimes x$. As $A^3\neq 0$, we have $x,y,z\neq 0$. Further, the uniqueness of a pure tensor decomposition for three non-zero factors implies there exists scalars $\alpha, \beta, \gamma$ with $$x=\alpha z,\quad y=\beta y,\quad z=\gamma x,\quad \alpha\beta\gamma=1.$$ As $\beta=1$, we get $\alpha\gamma=1$. As such, after translating the vectorized equalities back into square matrices, we find $A=\alpha A^2$ and $A^2=\gamma A$. It follows $$A^2\otimes A=(\gamma A)\otimes A=A\otimes (\gamma A)=A\otimes A^2.$$
\end{proof}

\begin{lemma}\label{lemma:square->scalar}
    Let $A$ be a nonzero $n\times n$ matrix, then $A^2\otimes A=A\otimes A^2$ if and only if $A^2=\mu A$ for $\mu$ a scalar.
\end{lemma}
\begin{proof}
    The converse was shown in the above proof, so we prove sufficiency. If $A^2\otimes A=A\otimes A^2$, then by Lemma \ref{lemma:scalar_prop}, assuming $A^2\neq 0$, there exists $\lambda\neq 0$ with $A=\lambda A^2$ and $A^2=\frac{1}{\lambda}A$. Letting $\mu=\frac{1}{\lambda}$, the result follows.
\end{proof}

\noindent \textbf{Proof of Theorem \ref{thm:R-Kronecker}.} By combining Lemma \ref{lemma:cubic->square} and Lemma \ref{lemma:square->scalar}, we have $$A^2\otimes A^3\otimes A=A\otimes A^3\otimes A^2\quad \iff \qquad A^2\otimes A=A\otimes A^2 \quad \iff \qquad A^2=\mu A$$ for $\mu$ a scalar.\\
\begin{flushright}\qed\end{flushright}



We now restrict our focus to adjacency matrices of a graph $\Gamma$. Equivalently, we can interpret these graph adjacency matrices as adjacency matrices of a quiver $Q$ with bidirectional arrows.

In this section, we assume $A^3\neq 0$. From Theorem \ref{thm:R-Kronecker}, it suffices to find graphs whose adjacency matrix satisfies $A^2=\mu A$ for some scalar $\mu$. We consider three different cases: The case where $\Gamma$ is without loops and multiple edges (Proposition \ref{prop:noloopsnomultiple}), the case where $\Gamma$ has loops but no multiple edges (Proposition \ref{prop:loopsnomultiple}) and lastly, when $\Gamma$ has loops and multiple edges (Theorem \ref{thm:2.9prime}).

\begin{proposition}\label{prop:noloopsnomultiple}
    Suppose $A$ is an adjacency $n\times n$ matrix corresponding to a graph $\Gamma$ without any loops or multiple edges, then $A^2=\mu A$ if and only if $A$ is the zero matrix, i.e. $\Gamma$ consists entirely of isolated vertices.
\end{proposition}
\begin{proof}
    As $\Gamma$ does not have any loops, we have $A_{ii}=0$ for all $i$, further each $A_{ij}$ is either $0$ or $1$. By Theorem \ref{thm:R-Kronecker}, it follows $$0=\mu A_{ii}=(A^2)_{ii}=\sum_{j=1}^nA_{ij}A_{ji}=\sum_{j=1}^nA_{ij}=\text{deg}(i).$$ We have thus shown $\deg(i)=0$ for all $i$, as such the graph $\Gamma$ has no edges and $A$ is the zero matrix.
\end{proof}

\begin{proposition}\label{prop:loopsnomultiple}
     Suppose $A$ is an adjacency $n\times n$ matrix corresponding to a graph $\Gamma$ without multiple edges (but it possibly has loops), then $A^2=\mu A$ if and only if $\Gamma$ is a disjoint union of complete graphs with $\mu$ vertices, where $\mu$ is the scalar satisfying $A^2=\mu A$.
\end{proposition}
\begin{proof}
    Suppose there exists a scalar $\mu$ such that $A^2=\mu A$. As $\Gamma$ does not have multiple edges, $A$ is a symmetric matrix with entries being either $0$ or $1$. Further, as we are allowing loops, we can have $A_{ii}\neq 0$. Similar to the proof of the previous proposition, we have $$\mu A_{ii}= (A^2)_{ii}=\sum_{j=1}^n A_{ij}A_{ji}=\sum^n_{j=1} A_{ij}=\deg(i)$$ where $\deg(i)$ is equal to the number of neighbors of $i$ (including loops). As $A^2=\mu A$, it follows $\deg(i)=0$, if $A_{ii}=0$, or $\deg(i)=\mu$, if $A_{ii}=1$. As such, we can partition our vertices of $A$ into a set of isolated vertices and the set $S=\{i:A_{ii}=1\}$, where every element $i$ of $S$ has $\deg(i)=\mu$.

    Restricting our attention to off-diagonal entries, we have $(A^2)_{ij}=\mu A_{ij}$, so if $A_{ij}=0$, then $A_{ij}^2=0$ and there are no $2$-paths from $i$ to $j$. If $A_{ij}=1$, then $(A^2)_{ij}=\mu$, so there are $\mu$ common neighbors of $i$ and $j$. As a corollary, the induced subgraph on $S$ is a disjoint union of complete graphs (with loops), each with $\mu$ vertices.
\end{proof}

\begin{theorem}\label{thm:2.9prime}
Let $K$ be a field of characteristic zero and let $A\in M_n(K)$ be a symmetric matrix with
$A^2=\mu A$ for some $\mu\in K$.
\begin{enumerate}
  \item If $\mu=0$, then $A^2=0$. If $K$ is formally real (in particular
        $K\subseteq\mathbb{R}$), this forces $A=0$, so all vertices are isolated.
  \item If $\mu\neq 0$, set $P=\mu^{-1}A$. Then $P$ is a symmetric idempotent,
        hence diagonalizable with all eigenvalues in $\{0,1\}$, and
        \[
            \operatorname{rank}A \;=\; \operatorname{rank}P
            \;=\; \operatorname{tr}P \;=\; \mu^{-1}\operatorname{tr}A .
        \]
        After a suitable reordering of vertices,
        \[
            A=\operatorname{diag}(A_1,\dots,A_m)\oplus 0,
        \]
        where each $A_i$ is the block on a connected component $C_i$ of the
        nonzero pattern of $A$, satisfies $A_i^2=\mu A_i$, and has
        \[
            \operatorname{rank}A_i \;=\; r_i \;=\; \mu^{-1}\operatorname{tr}A_i,
            \qquad 1\le r_i\le n_i=|C_i|.
        \]
        The block $A_i$ is of rank one \emph{(}i.e.\ $A_i=c_i\,w_iw_i^{\mathsf T}$
        for some $c_i\in K^{\times}$, $w_i\in K^{n_i}$\emph{)} if and only if
        $\operatorname{tr}A_i=\mu$.
\end{enumerate}
Conversely, for any symmetric idempotent $P\in M_n(K)$ the matrix $A:=\mu P$
satisfies $A^2=\mu A$.
\end{theorem}

\begin{proof}
\emph{(1)} If $\mu=0$ then $A^2=\mu A=0$. When $K$ is formally real,
$(A^2)_{ii}=\sum_{j}A_{ij}A_{ji}=\sum_j A_{ij}^2=0$ (using $A^{\mathsf T}=A$)
forces every $A_{ij}=0$, so $A=0$.

\emph{(2)} Since $A^2=\mu A$ and $\mu\ne 0$,
$P^2=\mu^{-2}A^2=\mu^{-2}(\mu A)=\mu^{-1}A=P$, and $P^{\mathsf T}=\mu^{-1}A^{\mathsf T}=P$,
so $P$ is a symmetric idempotent. Its minimal polynomial divides
$x(x-1)$, which has the distinct roots $0,1$; hence $P$ is diagonalizable with
eigenvalues in $\{0,1\}$ and $\operatorname{rank}P=\operatorname{tr}P$. As $\mu\ne0$,
$\operatorname{rank}A=\operatorname{rank}P$ and $\operatorname{tr}A=\mu\operatorname{tr}P$,
giving $\operatorname{rank}A=\mu^{-1}\operatorname{tr}A$.

Form the (undirected) graph on $\{1,\dots,n\}$ with an edge $\{i,j\}$ whenever
$A_{ij}\ne 0$. A symmetric matrix is block diagonal with respect to the
connected components of this graph; reordering vertices accordingly yields
$A=\operatorname{diag}(A_1,\dots,A_m)\oplus 0$, the final $0$ collecting the
isolated vertices (zero rows and columns). Each $A_i$ is symmetric, and because
the relation $A^2=\mu A$ respects the block decomposition we have $A_i^2=\mu A_i$;
applying the rank computation above to each block gives
$\operatorname{rank}A_i=\mu^{-1}\operatorname{tr}A_i=:r_i$. Finally $A_i$ has rank one
iff $r_i=1$ iff $\operatorname{tr}A_i=\mu$. In that case a rank-one symmetric matrix
is of the form $A_i=c_i\,w_iw_i^{\mathsf T}$, and $A_i^2=\mu A_i$ gives
$c_i\,(w_i^{\mathsf T}w_i)=\mu$ (in particular $w_i$ is non-isotropic).

The converse is immediate: $(\mu P)^2=\mu^2P^2=\mu^2P=\mu(\mu P)$.
\end{proof}

\begin{proposition}\label{prop:2.10prime}
Let $A\in M_n(\mathbb{R})$ be symmetric with \emph{nonnegative} entries and
$A^2=\mu A$.
\begin{enumerate}
  \item If $\mu=0$ then $A=0$.
  \item If $\mu>0$, then after reordering vertices
        \[
            A=\operatorname{diag}(A_1,\dots,A_m)\oplus 0,
            \qquad A_i=a_i a_i^{\mathsf T},
        \]
        where each $a_i\in\mathbb{R}_{>0}^{\,n_i}$ has strictly positive entries
        and $\lVert a_i\rVert^2=\mu$. Thus every nontrivial connected component of
        $Q$ is a rank-one complete weighted quiver. When $A$ is a $0/1$ matrix,
        each $a_i$ is the all-ones vector and $A_i=J_{n_i}$, recovering
        Proposition~\textup{2.8}: the components are complete graphs on $\mu$
        vertices with a loop at each vertex.
\end{enumerate}
\end{proposition}

\begin{proof}
Part (1) is Theorem~\ref{thm:2.9prime}(1). For (2), set $P:=\mu^{-1}A$, a
symmetric, entrywise nonnegative idempotent (hence an orthogonal projection).
Decompose into connected components as in Theorem~\ref{thm:2.9prime}: each
nonzero block $P_i$ is symmetric, nonnegative, idempotent, and \emph{irreducible}
(its pattern is connected). By the Perron--Frobenius theorem an irreducible
nonnegative matrix has spectral radius equal to a \emph{simple} eigenvalue with a
strictly positive eigenvector. Since $P_i$ is idempotent its eigenvalues lie in
$\{0,1\}$, and $P_i\ne 0$ forces the spectral radius to be $1$; simplicity then
makes $1$ an eigenvalue of multiplicity one and the only nonzero eigenvalue.
As $P_i$ is symmetric (so diagonalizable with $\operatorname{rank}=\operatorname{tr}$),
$\operatorname{rank}P_i=1$. Let $w_i$ be the positive unit Perron eigenvector;
then $P_i=w_iw_i^{\mathsf T}$. Putting $a_i:=\sqrt{\mu}\,w_i$ (possible since
$\mu>0$) gives $A_i=\mu P_i=a_i a_i^{\mathsf T}$ with $a_i\in\mathbb{R}_{>0}^{n_i}$ and
$\lVert a_i\rVert^2=\mu\lVert w_i\rVert^2=\mu$.

If moreover $A$ has $0/1$ entries, then $a_{i,k}a_{i,\ell}=(A_i)_{k\ell}\in\{0,1\}$;
since all $a_{i,k}>0$ and $a_{i,k}^2=(A_i)_{kk}\in\{0,1\}$ we get $a_{i,k}=1$ for
all $k$, so $A_i=J_{n_i}$ and $\mu=\lVert a_i\rVert^2=n_i$.
\end{proof}

\begin{example}

Consider the integer vector $a=\begin{pmatrix}
    1 \\2
\end{pmatrix}$ and the corresponding matrix $$A=a\cdot a^T=\begin{pmatrix}
    1\\2
\end{pmatrix}\begin{pmatrix}
    1 && 2
\end{pmatrix}=\begin{pmatrix}
    1 && 2\\
    2 && 4
\end{pmatrix}.$$ The above theorem guarantees $A$ satisfies the Quantum Yang-Baxter Equation with $A^2=\mu A$, where $\mu=||a||^2= 1^2+2^2=5$. Quickly checking, we see $A^2=\begin{pmatrix}
    5 && 10\\
    10 && 20
\end{pmatrix}=5A$, as desired. The theorem suggests all adjacency matrices satisfying the QYBE can be derived in a similar way. The corresponding graph for $A$ is pictured below:
\begin{figure*}[h]
    \centering
    \includegraphics[scale=.4]{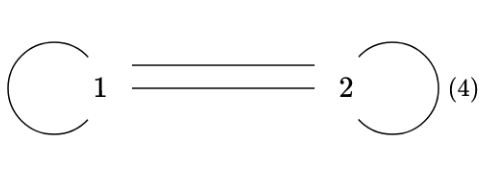}
\end{figure*}

\end{example}

\begin{corollary}\label{cor:2.11prime}
Let $Q$ be a finite quiver with symmetric adjacency matrix $A$ over a field $K$,
and assume $A^3\ne 0$. Then $Q$ satisfies the quantum Yang--Baxter equation if and
only if $A^2=\mu A$ for some $\mu\in K$, equivalently $A=\mu P$ with $P$ a symmetric
idempotent. In that case, after reordering vertices,
\[
    A=\operatorname{diag}(A_1,\dots,A_m)\oplus 0,
\]
each nonzero block $A_i$ corresponding to a connected component, with
$A_i^2=\mu A_i$ and $\operatorname{rank}A_i=\mu^{-1}\operatorname{tr}A_i$.

If, in addition, $A$ is a nonnegative real matrix \emph{(}e.g.\ a genuine
weighted graph\emph{)}, then by Proposition~\ref{prop:2.10prime} each nonzero
block is rank one, $A_i=a_ia_i^{\mathsf T}$ with $a_i\in\mathbb{R}_{>0}^{n_i}$ and
$\lVert a_i\rVert^2=\mu$; that is, every nontrivial connected component of $Q$ is a
complete weighted quiver. In general the blocks are symmetric idempotents
(scaled by $\mu$) of arbitrary rank, and the rank-one ``Temperley--Lieb'' case
occurs exactly for the components with $\operatorname{tr}A_i=\mu$.
\end{corollary}

\begin{proof}
By the computation of Section~2 and Theorem~2.3, the QYBE for $Q$ is equivalent
(under $A^3\ne0$) to $A^2=\mu A$. The structural statements are then
Theorem~\ref{thm:2.9prime}(2) and Proposition~\ref{prop:2.10prime}.
\end{proof}


\noindent Next, we proceed to show that a connected quiver with a groupoid structure satisfies the QYBE. Recall that a groupoid over a set $M$ is a set $G$ equipped with source and target mappings $\alpha, \beta: G \rightarrow M$, a multiplication map $m$ from $G_2=\{(g, h) \in G \times G: \beta(g)=\alpha(h)\}$ to $G$, an injective units mapping $\epsilon: M\rightarrow G$, and an inversion mapping $i:G\rightarrow G$, satisfying the following properties (where we write $gh$ for $m(g, h)$ and $g^{-1}$ for $i(g)$):

\begin{enumerate}
\item $g(hk)=(gh)k$.
\item $\epsilon(\alpha(g))g=g=g\epsilon(\beta(g))$.
\item $gg^{-1}=\epsilon(\alpha(g)$ and $g^{-1}g=\epsilon(\beta(g))$.
\end{enumerate}

\begin{proposition}
    Let $Q$ be a quiver arising from a groupoid with connected components $C_1, C_2,\ldots C_k$. Let $A$ be the adjacency matrix with $A_{C_1}, A_{C_2},\ldots, A_{C_k}$ the block matrices corresponding to the connected components. With this set-up, $A_{C_i}\otimes A_{C_i}$ satisfies the QYBE for each $i$. Equivalently, there exists a scalar $\mu_i$ such that $A^2_{C_i}=\mu_i A_{C_i}$. Furthermore, $A^2=\mu A$ for a single scalar $\mu$ if and only if $\mu_{1}=\mu_2=\ldots=\mu_k$. 
\end{proposition}

\begin{proof}
 Let $\Lambda$ be the object set of the groupoid $\mathcal{G}$ arising from $Q$. For our groupoid $\mathcal{G}$, let $$g_x:=|\text{Aut}(x)|\qquad (x\in \Lambda)$$ 
 If $i$ and $j$ lie in the same connected component, then $\text{Hom}(i,j)$ is nonempty and $$|\text{Hom}(i,j)|=|\text{Aut}(j)|$$ thus $$A_{ij}=|\text{Hom}(i,j)|=g_{j}$$ 
 If $i$ and $j$ are in different components, then $A_{ij}=0$, so each connected component with objects indexed by $C_i:=\{i_1,i_2,\ldots,i_n\}$ produces a block $$A_{C_i}=\textbf{1}_n\cdot (g_{i_1},\ldots,g_{i_n}),$$ where $\textbf{1}_n$ is the $n\times 1$ column vector consisting of ones. Thus, for any connected component $C$ the corresponding block $A_{C}$ is rank one with $A_{C}=\textbf{1}\cdot g^T$. 
Lastly, $$A_{C_i}^2=(\textbf{1}\cdot g^T)(\textbf{1}\cdot g^T)=\textbf{1}(g^T\cdot \textbf{1}) g^T=(\sum_{x\in C_i} g_x) \textbf{1}\cdot g^T=\mu_iA_{C_i}.$$ As such, each connected component satisfies $A^2\otimes A^3\otimes A=A\otimes A^3\otimes A^2$ by Theorem \ref{thm:R-Kronecker}. The observations regarding the global matrix $A$ follows immediately from the above observations.

\end{proof}

\begin{remark}
    As illustrated in the proof of the previous proposition, the scalar $\mu_i$ with $A_{C_i}^2=\mu_i A_{C_i}$ is the sum of the automorphism-group sizes in the $C_i$ connected component.  
\end{remark}

\begin{corollary}
Let $Q$ be a connected quiver with a groupoid structure, then $Q$ satisfies the quantum Yang-Baxter equation. 
\end{corollary}

From the previous corollary, it follows every connected groupoid induces a bialgebra structure via the FRT construction.

\begin{example}
    Consider a quiver with two connected components with a single object in both and $g_1=1$ and $g_2=2$. The corresponding matrix is $$A=\begin{bmatrix}1 &&0\\0&&2\end{bmatrix}.$$ There is not a scalar $\mu$ such that $A^2=\mu A$, despite both connected components satisfying $A_{C_i}^2=\mu_i A$ for some scalar $\mu_i$. 
\end{example}

\bigskip

\section{Construction of Hecke $R$-matrices}

\noindent It is known that one can produce constant solutions to the quantum Yang-Baxter equation through the Temperley-Lieb algebra $TL_N$ \cite{Baxter, Temperley-Lieb}. The Temperley-Lieb algebra $(TL_N)$ has $N-1$ generators $\{X_1,\ldots, X_{N-1}\}$ with relations 
\begin{align*}
    &X_k^2=dX_k\\
    &X_kX_jX_k=X_k \qquad\qquad |j-k|=1\\
    &X_jX_k=X_kX_j\qquad\qquad |j-k|>1.
\end{align*}

\noindent Temperley-Lieb algebra played a central role in Vaughan Jones’ discovery of his new polynomial invariant of knots which triggered a plethora of mathematical developments relating knot theory, topological quantum field theory, and statistical physics. 

Here we follow Kulish's rank-one tensor construction \cite{Kulish} to produce a Temperley–Lieb generator $X$, and we determine exactly when the associated operator $qI\pm X$ is a Hecke $R$-matrix. As we show in Proposition 3.5, this happens precisely under a normalization constraint relating $q$ to the invariant $\mu(b)=\tr(b^Tb^{-1})$. 

Before describing $X$ in terms of quivers, we first remind the reader of the aforementioned construction of $X$ introduced by Kulish \cite{Kulish}. The recipe is as follows:

Fix a local realization of the $TL_N$ algebra $\{1, X_1, X_2, \ldots, X_{N-1}\}$ in the space $$\mathcal{H}_N=\bigotimes_{k=1}^N \mathbb{C}^n_k. $$
\begin{enumerate}
    \item To define the above realization, take an $n\times n$ invertible matrix $b$, with inverse $\overline{b}=b^{-1}$.
    \item Through vectorization, view $b, \overline{b}$ as vectors of the $n^2$ dimensional space $\mathbb{C}^n\otimes \mathbb{C}^n$ with components $b_{ij}, i,j=1,2,\ldots, n$.
    \item Define the rank-one, Temperley-Lieb generator $X\in \text{End}(\mathbb{C}^n\otimes \mathbb{C}^n)$ as $$\text{vec}(b)\otimes \text{vec}(\overline{b}),\qquad X_{ab;cd}= b_{ab}\overline{b}_{cd}.$$ $X$ is a generator of $TL_N$ with $X_k$ acting on $\mathbb{C}^{n}_k\otimes \mathbb{C}^{n}_{k+1}$ of $\mathcal{H}_N$.
    \item $\Rc = qI + X$ satisfies the braid QYBE and the Hecke condition if and only if
$\mu(b) = \tr(b^{\mathsf T}b^{-1}) = -(q+q^{-1})$; more generally $\Rc = qI \pm X$ works under
$\mu(b) = \mp(q+q^{-1})$ as we will see below. 
\end{enumerate}

Let $\Bbbk$ be a field, $n \ge 2$, $V = \Bbbk^n$ with basis
$e_1,\dots,e_n$, and $q \in \Bbbk^\times$ with $q \neq \pm 1$. We use \emph{row-major} vectorization: for $b = (b_{ij}) \in M_n(\Bbbk)$ put
\[
  \operatorname{vec}(b) \;=\; \sum_{i,j=1}^{n} b_{ij}\, e_i \otimes e_j
  \;\in\; V \otimes V .
\]

\begin{definition}\label{def:Xb}
For $b \in GL_n(\kk)$ with $\bar b := b^{-1}$, set $\psi := \vect(b)$, $\varphi := \vect(\bar b)$ and
\[
  X = X(b) =\; \psi\,\varphi^{\mathsf T} \;\in\; \End(V \otimes V),
  \qquad X_{(i,j),(k,l)} = b_{ij}\,\bar b_{kl},
\]
i.e.\ $X(e_k \otimes e_l) = \bar b_{kl}\, \psi$. Put
\[
 \mu(b) =\; \sum_{i,j} b_{ij}\bar b_{ij} \;=\; \tr\!\left(b^{\mathsf T} b^{-1}\right).
\]
\end{definition}

Note $X(b)$ has rank one, $\tr X = \mu(b)$, and $X(\lambda b) = X(b)$ for $\lambda \in \kk^\times$,
so $\mu$ is a scale-invariant of $b$. Also $\mu(b^{\mathsf T}) = \mu(b)$.

\begin{lemma} \label{lem:TL}
Let $b \in GL_n(\kk)$ and $X = X(b)$. Write $X_1 = X \otimes I_V$ and $X_2 = I_V \otimes X$ on
$V^{\otimes 3}$. Then
\[
  X^2 = \mu(b)\,X, \qquad X_1 X_2 X_1 = X_1, \qquad X_2 X_1 X_2 = X_2 .
\]
Thus $X$ generates a Temperley--Lieb algebra with loop parameter $d = \mu(b)$, in the
normalisation $X_k^2 = d\,X_k$, $\;X_kX_jX_k = X_k$ for $|j-k| = 1$.
\end{lemma}

\begin{proof}
For the first relation,
$X^2(e_k \otimes e_l) = \bar b_{kl}\,X(\psi) = \bar b_{kl}\sum_{i,j} b_{ij}\,X(e_i \otimes e_j)
= \bar b_{kl}\Big(\sum_{i,j} b_{ij}\bar b_{ij}\Big)\psi = \mu(b)\,X(e_k \otimes e_l)$.

For the second, apply $X_1X_2X_1$ to $e_a \otimes e_b \otimes e_c$:
\begin{align*}
  X_1(e_a \otimes e_b \otimes e_c)
    &= \bar b_{ab}\sum_{i,j} b_{ij}\, e_i \otimes e_j \otimes e_c, \\
  X_2X_1(e_a \otimes e_b \otimes e_c)
    &= \bar b_{ab}\sum_{i,j} b_{ij}\,\bar b_{jc}\sum_{k,l} b_{kl}\, e_i \otimes e_k \otimes e_l
     = \bar b_{ab}\sum_{k,l} b_{kl}\, e_c \otimes e_k \otimes e_l,
\end{align*}
where the last step used $\sum_j b_{ij}\bar b_{jc} = (b\,b^{-1})_{ic} = \delta_{ic}$. Applying
$X_1$ once more and using $\sum_k \bar b_{ck} b_{kl} = (b^{-1}b)_{cl} = \delta_{cl}$ gives
$\bar b_{ab}\,\psi \otimes e_c = X_1(e_a \otimes e_b \otimes e_c)$. The relation
$X_2X_1X_2 = X_2$ is proved symmetrically, using $\bar b = b^{-1}$ on the other side.
\end{proof}

\begin{lemma}\label{lem:indep}
For $n \ge 2$ and $b \in GL_n(\Bbbk)$ the operators $X_1, X_2$ on $V^{\otimes 3}$ are distinct.
\end{lemma}

\begin{proof}
$\operatorname{Im}X_1 = \psi \otimes V$ and $\operatorname{Im}X_2 = V \otimes \psi$. Under the
flattening $V^{\otimes 3} \cong V \otimes (V \otimes V)$, an element $y \otimes \psi$ has rank $1$,
while $\psi \otimes x = \sum_i e_i \otimes \big(\sum_j b_{ij} e_j \otimes x\big)$ has rank
$\rank(b) = n \ge 2$ for $x \neq 0$. Hence the two images differ.
\end{proof}

\begin{proposition}\label{prop:normalisation}
Let $b \in GL_n(\kk)$, $n \ge 2$, $X = X(b)$, $\mu = \mu(b)$, and let $c \in \kk^\times$. Put
$\Rc := qI_{V \otimes V} + cX$. Then:
\begin{enumerate}
  \item $\Rc$ satisfies the braid form of the QYBE, $\Rc_{12}\Rc_{23}\Rc_{12} = \Rc_{23}\Rc_{12}\Rc_{23}$,
        if and only if
        \[
          q^2 + c\mu\,q + c^2 = 0 .
        \]
  \item $\Rc$ satisfies the Hecke condition $(\Rc - qI)(\Rc + q^{-1}I) = 0$ if and only if
        \[
          c\,\mu = -(q + q^{-1}).
        \]
  \item $\Rc$ satisfies both if and only if
        \[
          c = \pm 1 \qquad\text{and}\qquad \mu(b) = -c\,(q + q^{-1}).
        \]
\end{enumerate}
\end{proposition}

\begin{proof}
Write $Y_i = cX_i$, so that by Lemma~\ref{lem:TL} $Y_i^2 = \delta Y_i$ with $\delta := c\mu$, and
$Y_1Y_2Y_1 = c^3 X_1 = c^2 Y_1$, $Y_2Y_1Y_2 = c^2 Y_2$. Expanding,
\[
  \Rc_{12}\Rc_{23}\Rc_{12} = q^3 I + q^2(2Y_1 + Y_2) + q\,(Y_1Y_2 + Y_2Y_1) + q\delta Y_1 + c^2 Y_1,
\]
and the same expression with $Y_1 \leftrightarrow Y_2$ for $\Rc_{23}\Rc_{12}\Rc_{23}$. Subtracting,
\[
  \Rc_{12}\Rc_{23}\Rc_{12} - \Rc_{23}\Rc_{12}\Rc_{23} = \big(q^2 + q\delta + c^2\big)(Y_1 - Y_2).
\]
By Lemma~\ref{lem:indep}, $Y_1 \neq Y_2$, giving (1). For (2),
\[
  (\Rc - qI)(\Rc + q^{-1}I) = cX\big(cX + (q + q^{-1})I\big) = c\big(c\mu + q + q^{-1}\big)X,
\]
and $X \neq 0$, giving (2). For (3): substituting $c\mu = -(q+q^{-1})$ into $q^2 + c\mu q + c^2 = 0$
gives $q^2 - (q^2 + 1) + c^2 = 0$, i.e.\ $c^2 = 1$. Conversely $c = \pm1$ and
$c\mu = -(q+q^{-1})$ give both identities.
\end{proof}

\begin{corollary}\label{cor:nonempty}
Assume $\kk$ is algebraically closed and $n \ge 2$. Then $\mu : GL_n(\kk) \to \kk$ is surjective;
in particular, for every $q \in \kk^\times$ there exists $b \in GL_n(\kk)$ with
$\mu(b) = -(q+q^{-1})$, so that $\Rc = qI + X(b)$ is a Hecke $R$-matrix satisfying the braid QYBE.
\end{corollary}

\begin{proof}
Given a target value $\nu \in \kk$, choose $t \in \kk^\times$ with
$t + t^{-1} = -\nu + (n-2)$; such $t$ exists since $\kk$ is algebraically closed and the quadratic
$t^2 + (\nu - n + 2)t + 1 = 0$ has nonzero roots. Put
$b := \operatorname{diag}\big(b_t,\, I_{n-2}\big)$ with
$b_t = \left(\begin{smallmatrix} 0 & -t \\ 1 & 0\end{smallmatrix}\right)$. Then
$b^{-1} = \operatorname{diag}(b_t^{-1}, I_{n-2})$ and
\[
  \mu(b) = \tr\big(b_t^{\mathsf T}b_t^{-1}\big) + (n-2) = -(t + t^{-1}) + (n-2) = \nu .
\]
Taking $\nu = -(q+q^{-1})$ and applying Proposition~\ref{prop:normalisation}(3) with $c = 1$ establishes the claim. 
\end{proof}

\bigskip

\subsection*{Quiver Realizations of $X$}
From the above construction, matrices of the form $X=\text{vec}(b)\otimes \text{vec}(\overline{b})$ lead to nontrivial $R$-Hecke matrices. Following in the vein of the adjacency matrices of the previous section, we ask which quivers have adjacency matrices that can be expressed as $\text{vec}(b)\otimes \text{vec}(\overline{b})$.

To begin, $b$ is an arbitrary $n\times n$ invertible matrix, as such, $\text{vec}(b)\, \otimes\, \text{vec}(\overline{b})$ can have negative entries. With this in mind, it is more natural to consider weighted quivers, otherwise, the restrictions we impose on $b$ will be too rigid to lead to interesting Hecke $R$-matrices. Regardless, the weighted quivers we obtain have an underlying quiver structure that is similar to the quivers we have looked at in previous section; the connected components of the quiver are complete. 

The following theorem is a simple yet quite useful observation, but before stating the theorem, we introduce some terminology. We will call $b_{ij}$ the source entries and $\overline{b}_{kl}$ the target entries for reasons that will become apparent in the proof of the following theorem. Furthermore, as $b$ can be viewed as a weighted matrix, we may view any $b_{ij}$ as a weight $\alpha_{ij}$. We call such an $\alpha_{ij}$ a source weight. Similarly, we may view $\overline{b}_{kl}$ as a weight $\beta_{kl}$, which we call a target weight.
\begin{theorem} \label{recipe}
    Let $X$ be a $n^2\times n^2$ matrix. Then $X=\text{vec}(b)\,\otimes\,\text{vec}(\overline{b})$ for $b$ an invertible $n\times n$ matrix if and only if $X$ is a weighted adjacency matrix of a quiver $Q_X$ with the following structure:  
    \begin{itemize}
        \item $Q_X$ has $n^2$ vertices indexed as pairs $\{(i,j):1\leq i,j\leq n\}$.
       
        \medskip
        
        \item Edges exist between two vertices $(i, j)$ and $(k, l)$ exactly when the entries $b_{ij}$ and $\overline{b}_{kl}$ are both nonzero. Thus $Q_X$ is the complete directed quiver from the set of nonzero source-entry vertices $\{(i, j): b_{ij} \neq 0\}$ to the set of nonzero target-entry vertices $\{(k, l): \overline{b_{kl}} \neq 0\}$; these two sets need not coincide, and $Q_X$ is a complete graph precisely when they do.
        
        \medskip
        
        \item The weights of $Q_X$ factor multiplicatively as $w(\text{source})\cdot w(\text{target})$. Additionally, the source weights and the target weights form inverse matrices.
    \end{itemize}

    \noindent Lastly, we observe the following four-cycle equation holds: $$X_{(i,j),(k,l)}\cdot X_{(i',j'),(k',l')}=X_{(i,j),(k',l')}\cdot X_{(i',j'),(k,l)}$$
\end{theorem}

\begin{proof}
    First, we suppose $X=vec (b)\,\otimes vec(\overline{b})$ for $b$ an invertible $n\times n$ matrix. We write $u=vec(b)$ and $v=vec(\overline{b})$, thus $X=uv^T$, where $u,v$ are column vectors with $n^2$ entries and $X$ is an $n^2\times n^2$ matrix. 

    Regarding the first bullet point, suppose we have a basis $\{e_1, e_2,\ldots, e_n\}$ for $\mathbb{C}^n$, then $\{e_i\otimes e_j\}_{i,j=1}^n$ is a basis of $\mathbb{C}^n\otimes \mathbb{C}^n$, meaning the basis vectors are indexed by ordered pairs $(i,j)$. In particular, any operator $X:\mathbb{C}^n\otimes \mathbb{C}^n\rightarrow \mathbb{C}^n\otimes \mathbb{C}^n$ is represented by an $n^2\times n^2$ matrix whose rows and columns are indexed by pairs $(i,j)$. In the $\{e_i\otimes e_j\}_{i,j=1}^n$ basis, we may write $vec(b)=\displaystyle{\sum_{i,j}^n b_{i,j}e_i\otimes e_j}$ and $ X_{(i,j), (k,l)}=b_{ij}\cdot \overline{b}_{kl}$. As the basis vectors of $\mathbb{C}^n\otimes \mathbb{C}^n$ are $e_i\otimes e_j$, the vertices of $Q_X$ can be expressed as $$\text{vertices}(Q_X)=\{(i,j)\,\mid\, 1\leq i,j\leq n\}.$$

    For the second bullet point, the entries of $X_{(i,j), (k,l)}$ are defined as $$X_{(i,j), (k,l)}=b_{ij}\overline{b}_{kl};$$ therefore, if $b_{ij}=0$, then $X_{(i,j), (k,l)}=0$ for all $(k,l)$, meaning vertex $(i,j)$ has not outgoing arrows. Similarly, if $\overline{b}_{kl}=0$, then $(k,l)$ has no incoming arrows. This is the intuition behind calling $b_{ij}$ source entry vertices and $\overline{b}_{kl}$ target entry vertices. Thus in order to have an arrow $(i,j)\rightarrow (k,l)$, one must have both $b_{ij}\neq 0$ and $\overline{b}_{kl}\neq 0$.

    The observation that $Q_X$ is the complete directed quiver from the set of nonzero source-entry vertices $\{(i, j): b_{ij} \neq 0\}$ to the set of nonzero target-entry vertices $\{(k, l): \overline{b_{kl}} \neq 0\}$ follows immediately from the above realization.

    For the final bullet point, suppose we have an arrow $(i,j)\rightarrow (k,l)$, then if we label $(i,j)$ with the scalar $\alpha_{ij}=u(i,j)=b_{ij}$ and the vertex $(k,l)$ with the scalar $\beta_{k,l}=v(k,l)=\overline{b}_{kl}$, then the weight of the arrow $(i,j)\rightarrow (k,l)$ is $$w((i,j),(k,l))=X_{(i,j), (k,l)}=b_{ij}\cdot \overline{b}_{kl}=\alpha_{ij}\cdot \beta_{k,l}=w((i,j))\cdot w((k,l)).$$ By definition, we have $b=(\alpha_{ij})$ and $\overline{b}=(\beta_{kl})$, meaning the source weights and the target weights form inverse matrices.

    The converse of the theorem follows by defining $b=(\alpha_{ij})$ and $\overline{b}=(\beta_{kl})$; the structure is entirely determined by the weights.

    The four-cycle identity is immediate from definition:
    \begin{align*}
        X_{(i,j),(k,l)}\cdot X_{(i',j'), (k',l')} &= b_{ij}\cdot\overline{b}_{kl}\cdot b_{i'j'}\cdot \overline{b}_{k'l'}\\
        &=b_{ij}\cdot \overline{b}_{k'l'}\cdot b_{i'j'}\cdot\overline{b}_{kl}\\
        &=X_{(i,j),(k',l')}\cdot X_{(i',j'), (k,l)}.
    \end{align*} 
\end{proof}

\begin{example} \rm
    Consider the following weighted quiver:

    \begin{figure*}[h]
        \centering
        \includegraphics[scale=.35]{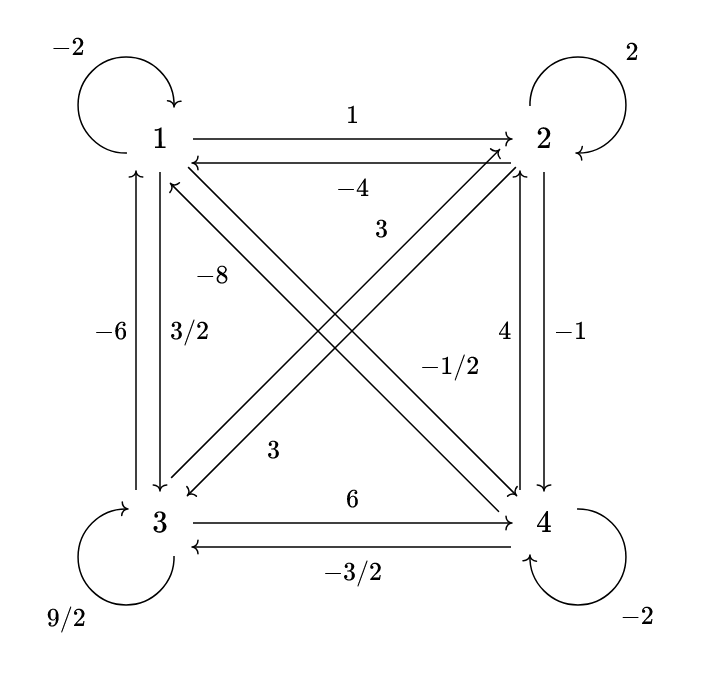}
    \end{figure*}
    
    \medskip

    \noindent The weighted adjacency matrix is given by $$X=\begin{pmatrix}
        -2 && 1 && 3/2 && -1/2\\
        -4 && 2 && 3 && -1\\
        -6 && 3 && 9/2 && -3/2\\
        -8 && 4 && 6 && -2
    \end{pmatrix}$$ which can be expressed as $$X = \text{vec}(b)\otimes \text{vec}(\overline{b})$$ where $$b=\begin{pmatrix}
        1 && 2\\
        3 && 4
    \end{pmatrix}$$ and $$\overline{b}=
        -1/2\begin{pmatrix}
            4 && -2\\
            -3 && 1
        \end{pmatrix}$$

\medskip

Here
\[
  \mu(b) = \tr(b^{\mathsf T}b^{-1}) = \tr X = -2 + 2 + \tfrac92 - 2 = \tfrac52 ,
\]
so $X^2 = \tfrac52 X$ and $X$ is a Temperley--Lieb generator with $d = 5/2$. 

\bigskip

By Proposition~\ref{prop:normalisation}, $qI + X$ is a Hecke $R$-matrix precisely for
$q + q^{-1} = -\tfrac52$, i.e.\ $q \in \{-2, -\tfrac12\}$, and $qI - X$ precisely for
$q \in \{2, \tfrac12\}$. 

For all other $q$ the operator $qI \pm X$ satisfies neither the Hecke
condition nor the braid relation. The four-cycle identity
$X_{(i,j),(k,l)}X_{(i',j'),(k',l')} = X_{(i,j),(k',l')}X_{(i',j'),(k,l)}$ holds unconditionally,
being equivalent to $\rank X = 1$.

 The main point is that the product of the weights of the red arrows is equal to the product of the weights of the blue arrows.

    \begin{figure*}[h]
        \centering
        \includegraphics[scale=.35]{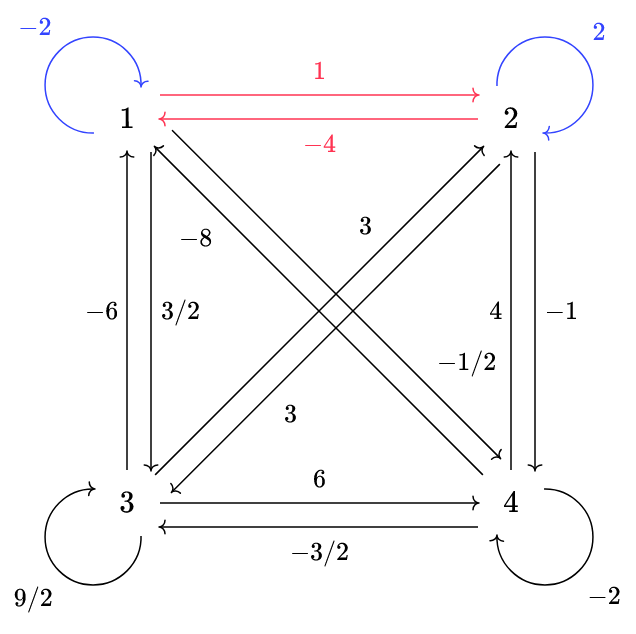}
    \end{figure*}
    
\end{example}

\begin{example} \label{Toeplitz} \rm

\noindent Let $Q$ be the Toeplitz quiver below:

\begin{figure*}[h]
    \includegraphics[scale=.35]{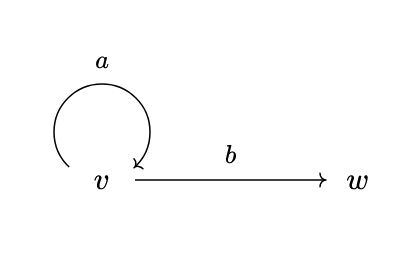}
\end{figure*}

\bigskip

Let $V = \operatorname{span}\{a,b\}$. The classical antisymmetrizer
$\psi_{\mathrm{cl}} := a \otimes b - b \otimes a$ equals
$\operatorname{vec}(b_{\mathrm{cl}})$ for
\[
  b_{\mathrm{cl}} = \begin{pmatrix} 0 & 1 \\ -1 & 0 \end{pmatrix},
  \qquad
  \bar b_{\mathrm{cl}} = b_{\mathrm{cl}}^{-1}
    = \begin{pmatrix} 0 & -1 \\ 1 & 0 \end{pmatrix},
\]
whence
\[
  X_{\mathrm{cl}}
    = \operatorname{vec}(b_{\mathrm{cl}})\operatorname{vec}(\bar b_{\mathrm{cl}})^{\mathsf T}
    = \begin{pmatrix}
        0 & 0 & 0 & 0 \\ 0 & -1 & 1 & 0 \\ 0 & 1 & -1 & 0 \\ 0 & 0 & 0 & 0
      \end{pmatrix},
  \qquad
  \mu(b_{\mathrm{cl}}) = -2 .
\]
In particular $X_{\mathrm{cl}}^2 = -2X_{\mathrm{cl}}$. By Proposition 3.5,
$qI + X_{\mathrm{cl}}$ satisfies the Hecke condition only for $q + q^{-1} = 2$,
i.e.\ $q = 1$, and $qI - X_{\mathrm{cl}}$ only for $q = -1$; both are excluded by
the standing hypothesis $q \neq \pm 1$. Thus the undeformed antisymmetrizer
produces no Hecke $R$-matrix, and a $q$-deformation is forced.

Accordingly, set
\[
  \psi_q := a \otimes b - q\, b \otimes a = \operatorname{vec}(b_q),
  \qquad
  b_q = \begin{pmatrix} 0 & -q \\ 1 & 0 \end{pmatrix},
  \qquad
  \bar b_q = \begin{pmatrix} 0 & 1 \\ -q^{-1} & 0 \end{pmatrix},
\]
for which
\[
  X_q = \operatorname{vec}(b_q)\operatorname{vec}(\bar b_q)^{\mathsf T}
      = \begin{pmatrix}
          0 & 0 & 0 & 0 \\ 0 & -q & 1 & 0 \\ 0 & 1 & -q^{-1} & 0 \\ 0 & 0 & 0 & 0
        \end{pmatrix},
  \qquad \mu(b_q) = -(q + q^{-1}).
\]
Hence $\check R_2(q) = qI + X_q$ satisfies the braid QYBE and the Hecke condition
for every $q \in k^\times$, and equals the standard $GL_q(2)$ check-$R$-matrix
\[
  \check R_2(q) = \begin{pmatrix}
    q & 0 & 0 & 0 \\ 0 & 0 & 1 & 0 \\ 0 & 1 & q - q^{-1} & 0 \\ 0 & 0 & 0 & q
  \end{pmatrix},
\]
recovering Lemma 3.11. Setting $q = 1$ returns $X_{\mathrm{cl}}$ and the classical
flip.

\end{example}

\subsection*{The Temperley--Lieb projector computation}

Let $X = X_{\mathrm{cl}}$ be as in Example~\ref{Toeplitz}, so that $X^2 = -2X$. Then
\[
  P := -\tfrac12 X \quad\text{satisfies}\quad P^2 = \tfrac14 X^2 = -\tfrac12 X = P,
\]
so $P$ is an idempotent. Set $e := \alpha P$ with $\alpha \in \kk^\times$;
we test the Temperley--Lieb relations
\[
  \text{(i) } e^2 = \alpha e, \qquad
  \text{(ii) } e_ie_je_i = e_i \ \ (|i-j| = 1), \qquad
  \text{(iii) } e_ie_j = e_je_i \ \ (|i-j| > 1).
\]
Relation (i) holds for all $\alpha$ since $e^2 = \alpha^2P^2 = \alpha^2 P = \alpha e$, and (iii) is
automatic. For (ii), Lemma~\ref{lem:TL} gives $X_1X_2X_1 = X_1$ and therefore
\[
  P_{12}P_{23}P_{12} = \left(-\tfrac12\right)^3 X_1X_2X_1 = -\tfrac18 X_1 = \tfrac14 P_{12}.
\]
Hence
$e_{12}e_{23}e_{12} = \alpha^3 P_{12}P_{23}P_{12} = \tfrac{\alpha^3}{4}P_{12}$, and (ii) requires
$\alpha P_{12} = \tfrac{\alpha^3}{4}P_{12}$, i.e.
\[
  \alpha^2 = 4, \qquad \alpha = \pm 2 .
\]
With the normalization $e = -(q + q^{-1})X = 2(q+q^{-1})P$ we have $\alpha = 2(q + q^{-1})$, so
\[
  \pm 2 = 2(q + q^{-1}) \iff q + q^{-1} = \pm 1 \iff q^2 \mp q + 1 = 0,
\]
with solutions
\[
  q = \frac{1 \pm i\sqrt{3}}{2} = e^{\pm i\pi/3}
  \qquad\text{and}\qquad
  q = \frac{-1 \pm i\sqrt{3}}{2} = e^{\pm 2i\pi/3}.
\]

\begin{remark}
Here $X$ is a rank-one Temperley--Lieb generator acting on $V \otimes V$ with
$V = \Bbbk^n$. A natural question is whether the standard $GL_q(N)$ $R$-matrix
(to be used in the next section) can be produced this way.

For $N = 2$ the answer is yes: the standard $GL_q(2)$ $R$-matrix is exactly
$qI + X$ with $X$ rank one, and it arises from Theorem 3.7 with
$b = \begin{pmatrix} 0 & -q \\ 1 & 0 \end{pmatrix}$.

For $N\ge 3$, a \emph{single} instance of Theorem \ref{recipe} cannot produce the
standard $GL_q(N)$ R-matrix, since for the standard R-matrix the operator
$X=\check R-qI$ has rank $\binom{N}{2}$, which exceeds one. However the
recipe does produce $\check R$ as a \emph{sum} of $\binom{N}{2}$ rank-one
pieces, one per unordered pair of basis vectors. We describe and prove this
realization below.
\end{remark}

\subsection*{The basic 2-dimensional building block}

\begin{lemma} \label{lem:N=2}
Let $V_2 = \Bbbk^2$ with basis $e_1, e_2$, and let
$b = \begin{pmatrix} 0 & -q \\ 1 & 0 \end{pmatrix} \in M_2(\Bbbk)$,
which is invertible with
$\bar b = b^{-1} = \begin{pmatrix} 0 & 1 \\ -q^{-1} & 0 \end{pmatrix}$. Set
\[
  X_0 \;=\; \operatorname{vec}(b)\,\operatorname{vec}(\bar b)^{\mathsf T}
  \;\in\; \operatorname{End}(V_2 \otimes V_2).
\]
Then $\check R_2(q) = q\,I_{V_2 \otimes V_2} + X_0$ is the standard $GL_q(2)$
Hecke $R$-matrix.
\end{lemma}

\begin{proof}
By the row-major convention, $\operatorname{vec}(b) = (0,-q,1,0)^{\mathsf T}$ and
$\operatorname{vec}(\bar b) = (0,1,-q^{-1},0)^{\mathsf T}$, both in the basis
$(e_1 \otimes e_1,\, e_1 \otimes e_2,\, e_2 \otimes e_1,\, e_2 \otimes e_2)$
of $V_2 \otimes V_2$. Hence
\[
  X_0 = \begin{pmatrix}
    0 & 0 & 0 & 0 \\ 0 & -q & 1 & 0 \\ 0 & 1 & -q^{-1} & 0 \\ 0 & 0 & 0 & 0
  \end{pmatrix},
  \qquad
  \check R_2(q) = qI + X_0 = \begin{pmatrix}
    q & 0 & 0 & 0 \\ 0 & 0 & 1 & 0 \\ 0 & 1 & q - q^{-1} & 0 \\ 0 & 0 & 0 & q
  \end{pmatrix},
\]
which acts on basis vectors as
\[
  \check R_2(q)(e_i \otimes e_j) =
  \begin{cases}
    q\, e_i \otimes e_i, & i = j, \\
    e_j \otimes e_i, & i < j, \\
    e_j \otimes e_i + (q - q^{-1})\, e_i \otimes e_j, & i > j,
  \end{cases}
\]
precisely the standard $GL_q(2)$ check-$R$-matrix. Moreover
$\mu(b) = \operatorname{tr}(b^{\mathsf T}b^{-1}) = -(q + q^{-1})$, so
$X_0^2 = -(q+q^{-1})X_0$ by Lemma 3.2, and $\check R_2(q)$ satisfies the braid
QYBE and the Hecke condition by Proposition 3.5 with $c = 1$.
\end{proof}

\subsection*{The standard $GL_q(N)$ R-matrix as a pairwise sum}

Let $V=\K^N$ with basis $e_1,\dots,e_N$, and for each unordered pair
$\{i,j\}\subset \{1,\dots,N\}$ with $i<j$, let
\[
V_{ij} = \Span_\K\{e_i,e_j\}\;\subset\; V
\]
be the 2-dimensional subspace spanned by $\{e_i,e_j\}$.

\begin{definition}\label{def:Xij}
For each pair $i<j$, define
\[
X^{ij}\;\in\;\End(V\otimes V)
\]
as follows. The subspace
\[
W_{ij} = \Span_\K\{e_i\otimes e_j,\ e_j\otimes e_i\}\;\subset\;V\otimes V
\]
is 2-dimensional. Choose the basis $(e_i\otimes e_j,\ e_j\otimes e_i)$ for
$W_{ij}$. Define $X^{ij}$ to act on $W_{ij}$ by the matrix
\[
[X^{ij}]_{W_{ij}} = \begin{pmatrix}-q & 1\\ 1 & -q^{-1}\end{pmatrix},
\]
and to vanish on every other basis vector $e_a\otimes e_b$ with
$\{a,b\}\ne\{i,j\}$.
\end{definition}

\begin{lemma}\label{lem:Xij-rank-one}
Each $X^{ij}$ is a rank-one Temperley--Lieb generator with loop parameter 
$\mu = -(q+q^{-1})$, and it is exactly of the form prescribed by
Theorem \ref{recipe} applied to the matrix
$b^{ij} = \left(\begin{smallmatrix}0&-q\\1&0\end{smallmatrix}\right)$
acting in the $(e_i,e_j)$ basis.
\end{lemma}

\begin{proof}
The $2\times 2$ matrix
$M=\begin{pmatrix}-q & 1\\ 1 & -q^{-1}\end{pmatrix}$
has $\det M = (-q)(-q^{-1})-1\cdot 1 = 1-1 = 0$, so $\operatorname{rank}M=1$.
Trace is $\mathrm{tr}\,M = -(q+q^{-1})$, so the nonzero eigenvalue is
$-(q+q^{-1})$. Hence $M^2=-(q+q^{-1})M$, so $X^{ij}$ extends this to
$(X^{ij})^2=-(q+q^{-1})X^{ij}$ on all of $V\otimes V$ (vanishing outside
$W_{ij}$ in both factors).

Identification with the recipe of Theorem \ref{recipe} in the 2-dim subspace $V_{ij}$:
take $b^{ij} = \left(\begin{smallmatrix}0&-q\\1&0\end{smallmatrix}\right)$ and apply
Lemma~\ref{lem:N=2} to get the same matrix $M$ on $W_{ij}$.
\end{proof}

\begin{theorem}\label{thm:main}
The standard $GL_q(N)$ $\check R$-matrix is
\begin{equation}\label{eq:main}
\check R_N(q) =\; qI_{V\otimes V} \;+\; \sum_{1\le i<j\le N} X^{ij}.
\end{equation}
The matrix $\check R_N(q)$ satisfies the Hecke
condition $(\check R_N-qI)(\check R_N+q^{-1}I)=0$ and the quantum
Yang--Baxter equation
$\check R_{12}\check R_{23}\check R_{12}=\check R_{23}\check R_{12}\check R_{23}$
on $V^{\otimes 3}$.
\end{theorem}

\begin{proof}
The image of the basis vector $e_a\otimes e_b$ under $qI_{V\otimes V} \;+\; \sum_{1\le i<j\le N} X^{ij}$
\eqref{eq:main} is:
\begin{itemize}[leftmargin=2em]
\item If $a=b$: only the $qI$ term contributes ($X^{ij}$ vanishes on
diagonal basis vectors). So
$qI_{V\otimes V} \;+\; \sum_{1\le i<j\le N} X^{ij}(e_a\otimes e_a) = q\,e_a\otimes e_a$.
\item If $a\ne b$: let $i=\min(a,b),\ j=\max(a,b)$. Then $e_a\otimes e_b\in
W_{ij}$, and only $X^{ij}$ (among the summands) acts non-trivially on it.
\begin{itemize}
\item If $a<b$ (i.e.\ $a=i,\ b=j$): in the chosen basis
$(e_i\otimes e_j,\ e_j\otimes e_i)$ of $W_{ij}$, the vector
$e_a\otimes e_b=e_i\otimes e_j$ has coordinate $(1,0)^T$. Applying
$X^{ij}=\begin{pmatrix}-q&1\\1&-q^{-1}\end{pmatrix}$ gives
$(-q,1)^T = -q\,e_i\otimes e_j + e_j\otimes e_i$. Adding $q(e_i\otimes
e_j)$ from the $qI$ contribution: $0\cdot e_i\otimes e_j + e_j\otimes
e_i = e_b\otimes e_a$.
So $qI_{V\otimes V} \;+\; \sum_{1\le i<j\le N} X^{ij}(e_a\otimes e_b)=e_b\otimes e_a$ for $a<b$. 
\item If $a>b$ (i.e.\ $a=j,\ b=i$): the vector
$e_a\otimes e_b=e_j\otimes e_i$ has coordinate $(0,1)^T$ in the basis. Applying
$X^{ij}$ gives $(1,-q^{-1})^T = e_i\otimes e_j - q^{-1} e_j\otimes e_i$.
Adding $q(e_j\otimes e_i)$ from $qI$: $e_i\otimes e_j + (q-q^{-1})
e_j\otimes e_i$.
So $qI_{V\otimes V} \;+\; \sum_{1\le i<j\le N} X^{ij}(e_a\otimes e_b)=e_b\otimes e_a + (q-q^{-1})\,e_a\otimes
e_b$ for $a>b$. 
\end{itemize}
\end{itemize}
These match the defining action of the standard $GL_q(N)$ check-R-matrix. Hence \eqref{eq:main} holds.

$\check R_N(q)$ is the standard $GL_q(N)$
R-matrix, which is well-known to satisfy both conditions
\cite[\S 9]{Klimyk-Schmudgen}.
\end{proof}

The decomposition~\eqref{eq:main} admits a clean quiver-theoretic
interpretation.

\begin{remark} \label{cor:quiver-realization} \rm
For each pair $i<j$, let $Q^{ij}$ be the Toeplitz quiver on the two
vertices $\{i,j\}$ (one loop at each vertex, and two parallel directed
edges between them) decorated with the weights determined by
$b^{ij} = \left(\begin{smallmatrix}0&-q\\1&0\end{smallmatrix}\right)$ via Theorem \ref{recipe}. Equivalently, $Q^{ij}$ is the underlying weighted quiver of
Example \ref{Toeplitz} with $q$-deformed weights.

The standard $GL_q(N)$ check-R-matrix is the sum
\[
\check R_N(q) \;=\; qI + \sum_{1\le i<j\le N} X_{Q^{ij}},
\]
where each summand $X_{Q^{ij}}$ is the rank-one Temperley--Lieb operator
attached by the construction discussed before to the quiver
$Q^{ij}$, embedded in $V\otimes V$ via the inclusion $V_{ij}\hookrightarrow V$.

The full \emph{configuration} of weighted quivers
$\{Q^{ij}\}_{1\le i<j\le N}$ that produces the standard $GL_q(N)$
R-matrix may be viewed as a single weighted quiver $Q_N$
on $N$ vertices, where every pair $\{i,j\}$ carries its own copy of the
Toeplitz weighted structure. In adjacency-matrix language,
$Q_N$ has the $N\times N$ adjacency matrix
\[
A_{Q_N} = \begin{pmatrix}
0 & w & w & \cdots & w\\
w & 0 & w & \cdots & w\\
\vdots & & \ddots & & \vdots\\
w & w & w & \cdots & 0
\end{pmatrix}
\]
\end{remark}

\section{Hecke-deformed face algebras of quivers}

\noindent Throughout, $\Bbbk$ is a field and $q\in^{\times}$ with $q\neq\pm1$; we also assume $q$ is not a root of unity. All quivers are
finite. The purpose of this section is to study Hayashi's face algebra over quivers with QYBE and Hecke condition. 

\subsection{Hayashi's face algebra}
Let $Q=(Q_0, Q_1, s, r)$ be a finite quiver where $Q_0$ is the set of vertices, $Q_1$ is the set of arrows and both $s$, and $r$ are functions from $Q_1$ to $Q_0$ such that for any arrow $e$ from $u$ to $v$, $s(e)=u$ and $r(e)=v$. The \textit{path algebra} of $Q$ over a field $\Bbbk$, denoted $\Bbbk Q$ is the $\Bbbk$-algebra with $\Bbbk$-basis given by all the paths in $Q$ and ring structure determined by path concatenation. We will read paths of $Q$ from left-to-right. For any vertex $i\in Q_0$, we associate a trivial path $e_i$. The path algebra $\Bbbk Q$ is $\mathbb N$-graded by path length, where $( \Bbbk Q)_k=\Bbbk Q_k$, with $Q_k$ consisting of paths of length $k\in \mathbb N$. 

\begin{definition}[{\cite[Ex.~1.1]{Hay96}}]\label{def:face} \rm
Hayashi's face algebra $\mathfrak H(Q)$ for a finite quiver $Q$ over a field $\Bbbk$ is generated by elements $\{x_{a, b}\}$, for $a, b \in Q_k$ subject to the following relations:

(1) $x_{i, j} x_{u, v}=\delta_{i, u} \delta_{j, v} x_{i, j}$, for each $i, j, u, v\in Q_0$.

(2) $x_{s(p), s(q)} x_{p, q} =x_{p, q} =x_{p, q} x_{r(p), r(q)}$, for each $p, q \in Q_1$. 

(3) $x_{p, q} x_{p', q'} =\delta_{r(p), s(p')} \delta_{r(q), s(q')} x_{pp', qq'}$, for each $p, q, p', q' \in Q_1$. 

\bigskip

\noindent Clearly, $\mathfrak H(Q)$ is a unital $\mathbb K$-algebra. The idempotents $x_{i,j}$, $i,j \in Q_0$, are
orthogonal with $1_{\mathfrak H(Q)} = \sum_{i,j} x_{i,j}$, and $x_{a,b} \in x_{s(a),s(b)}\,\mathfrak H(Q)\,x_{r(a),r(b)}$. 
\end{definition}

Recall that a \emph{weak bialgebra} $H$ has compatible algebra and coalgebra structures with
the unit and counit axioms relaxed; it carries two canonical coideal subalgebras,
the source and target counital subalgebras $H_s,H_t$, and $H$ is an ordinary
bialgebra exactly when $H_s=H_t=\Bbbk$ \cite[\S2]{HWWW}.

\begin{proposition} \cite[Prop.~3.4]{CW}
\label{prop:HQ}
There is an isomorphism of $\mathbb N$-graded algebras
$\mathfrak H(Q)\cong \mathbb K\widehat Q$, sending $x_{a,b}\mapsto[a,b]$. Moreover, by
\cite[Thm.~4.17]{HWWW}, $\mathfrak H(Q)$ is an
$\mathbb N$-graded weak bialgebra coacting
universally on the path algebra $\mathbb KQ$ with
\[
\Delta(x_{a,b})=\sum_{c}x_{a,c}\otimes x_{c,b},\qquad
\varepsilon(x_{a,b})=\delta_{a,b},
\]
the sum over paths $c$ of the same length as $a,b$.

\end{proposition}

The ring-theoretic and homological properties of the Hayashi's face algebra are completely
determined by graph data of $Q$ via $\widehat Q$ \cite{CW}; in particular
$\operatorname{GKdim}\mathfrak H(Q)=2\operatorname{GKdim} \mathbb KQ-1$ when the latter
is finite. 

\subsection{Deformation of face algebra}

\noindent Let $Q=(Q_0, Q_1, s, r)$ be a finite quiver. For each $u\in Q_0$, let $L_Q(u)=\Span_{\Bbbk}\{l \in Q_1: s(l)=r(l)=u\}$ be the $\Bbbk$-vector space of loops at $u$. Consider an invertible
$\check R_u\in\mathrm{GL}\big(L_Q(u)\otimes L_Q(u)\big)$ for each $u\in Q_0$ satisfying the following conditions:\\
\begin{itemize}
 \item [(1)] (Braid form of Quantum Yang-Baxter Equation) \[
(\check R_u\otimes I)(I\otimes\check R_u)(\check R_u\otimes I)
=(I\otimes\check R_u)(\check R_u\otimes I)(I\otimes\check R_u)
\]
 \item [(2)] (Hecke Condition) $\check R_u^{2}=I+(q-q^{-1})\check R_u$, equivalently
$(\check R_u-qI)(\check R_u+q^{-1}I)=0$.
\end{itemize}

\noindent Let $d_u := \dim_\kk L_Q(u)$. The loops of $\widehat Q$ at $[u,u]$ are $L_Q(u)\otimes L_Q(u)$. We arrange the matrix of loop generators as $X^{(u)}=\sum_{a,b}E_{ab}\otimes x_{[\ell_a,\ell_b]}$  with $\{\ell_a\}$  a basis of $L_Q(u)$.  
We impose RTT relation with using the associated FRT matrix $R_u:=\tau\circ\check R_u$ (with $\tau$ being the tensor flip):
\[
R_u\,X^{(u)}_1X^{(u)}_2=X^{(u)}_2X^{(u)}_1\,R_u .
\]
where $X_1^{(u)}=X^{(u)}\otimes I$ and $X_2^{(u)}=I\otimes X^{(u)}$. 

\begin{definition}\label{def:Iloop}
Let $I_{\mathrm{loop}} \subseteq \mathfrak H(Q)$ be the two-sided ideal generated by the entries of the local
RTT relations $R_u X^{(u)}_1X^{(u)}_2 - X^{(u)}_2X^{(u)}_1R_u$, $u \in Q_0$. Set
\[
  \mathfrak H^{\mathrm{def}}(Q) \;:=\; \mathfrak H(Q)/I_{\mathrm{loop}},
\]
\emph{a quotient algebra of $\mathfrak H(Q)$}.
\end{definition}

The ideal $I_{\mathrm{loop}}$ is, in general, \emph{not} a coideal. We give below some examples to illustrate the above construction and discuss cases when $I_{\mathrm{loop}}$ is a biideal.

\begin{definition}\label{def:rose}
The \emph{rose with $n$ petals} $\mathcal R_n$ is the quiver with one vertex $v$ and
$n$ loops $\ell_1,\dots,\ell_n$. Then $\kk\mathcal{R}_n = \kk\langle \ell_1,\dots,\ell_n\rangle$ is free, and $\widehat{\mathcal R_n}$ has the single vertex $[v, v]$ and the $n^2$
loops $x_{ab}:=[\ell_a,\ell_b]$, $1\le a,b\le n$. By Proposition~\ref{prop:HQ},
\[
\mathfrak H(\mathcal R_n)\cong \Bbbk \widehat{R_n}=\Bbbk\langle x_{ab}:1\le a,b\le n\rangle
\]
Assemble the loop generators into the matrix $X=(x_{ab})_{1\le a,b\le n}$. The Hecke-deformed face algebra of $R_n$ is defined as the quotient of the free face algebra by the ideal $I_{loop}$ which is generated by the entries of the RTT relation.
\[
  \mathfrak H^{\mathrm{def}}(\mathcal{R}_n) \;=\; \mathfrak H{\mathcal{R}}_n \big/
  \big(\text{entries of } RX_1X_2 - X_2X_1R\big).
\]
where $X_1=X\otimes I_n,\ \ X_2=I_n\otimes X$.

\end{definition}

\begin{convention}\label{conv:Rcomp}
Let $\{E_{ij}\}_{1\le i,j\le n}$ denote the standard matrix units in $M_n(\Bbbk)$. Let
$R\in\operatorname{End}(\Bbbk^n\otimes \Bbbk^n)$ be the standard $GL_q(n)$ Hecke
$R$-matrix
\begin{equation}\label{eq:R-matrix-def}
  R 
  \;=\;
  q \sum_{i=1}^n E_{ii} \otimes E_{ii}
  \;+\;
  \sum_{\substack{i,j=1 \\ i\neq j}}^n E_{ii} \otimes E_{jj}
  \;+\;
  (q - q^{-1}) \sum_{\substack{i,j=1 \\ i<j}}^n E_{ij} \otimes E_{ji}.
\end{equation}

Components of $R \in \End(\kk^n \otimes \kk^n)$ are indexed by
$R(e_i \otimes e_j) = \sum_{a,b} R^{ab}_{ij}\,e_a \otimes e_b$, so that $R^{ab}_{ij}$ is the
$\big((a,b),(i,j)\big)$ entry of $R$, rows indexed by superscripts.
\end{convention}

\begin{lemma}\label{lem:Rcomp}
For $R$ as in Equation (\ref{eq:R-matrix-def}), the nonzero components are
\[
  R^{ii}_{ii} = q, \qquad
  R^{ij}_{ij} = 1 \ (i \neq j), \qquad
  R^{ji}_{ij} = q - q^{-1} \ (i > j),
\]
and no others.
\end{lemma}

\begin{proof}
The proof is straightforward. 
\end{proof}

\begin{proposition} \label{prop:relations}
Under Convention~\ref{conv:Rcomp}, the entrywise form of $R\,X_1X_2 = X_2X_1\,R$ is
\begin{equation}\label{eq:RTTentry}
  \sum_{a,b} R^{ij}_{ab}\; x_{ak}x_{b\ell} \;=\; \sum_{a,b} x_{jb}x_{ia}\; R^{ab}_{k\ell},
  \qquad 1 \le i,j,k,\ell \le n .
\end{equation}
For $R$ as in Equation (\ref{eq:R-matrix-def}), the relations \eqref{eq:RTTentry} are equivalent to the following, for
all $i < k$ and $j < \ell$:
\begin{align}
  x_{i\ell}\,x_{ij} &= q\;x_{ij}\,x_{i\ell}, \tag{R1} \\
  x_{kj}\,x_{ij} &= q\;x_{ij}\,x_{kj}, \tag{R2} \\
  x_{i\ell}\,x_{kj} &= x_{kj}\,x_{i\ell}, \tag{R3} \\
  x_{k\ell}\,x_{ij} - x_{ij}\,x_{k\ell} &= (q - q^{-1})\;x_{i\ell}\,x_{kj}. \tag{R4}
\end{align}
\end{proposition}

\begin{proof}
With $(X_1)^{ab}_{cd} = x_{ac}\delta_{bd}$ and $(X_2)^{ab}_{cd} = \delta_{ac}x_{bd}$ one has
$(X_1X_2)^{ab}_{cd} = x_{ac}x_{bd}$ and $(X_2X_1)^{ab}_{cd} = x_{bd}x_{ac}$, whence
$(RX_1X_2)^{ij}_{k\ell} = \sum_{a,b}R^{ij}_{ab}x_{ak}x_{b\ell}$ and
$(X_2X_1R)^{ij}_{k\ell} = \sum_{a,b}x_{jb}x_{ia}R^{ab}_{k\ell}$, which is \eqref{eq:RTTentry}.

By Lemma~\ref{lem:Rcomp}, for fixed $(i,j)$ the only $(a,b)$ with $R^{ij}_{ab} \neq 0$ are
$(a,b) = (i,j)$, with coefficient $q$ if $i = j$ and $1$ otherwise, and $(a,b) = (j,i)$ when
$i < j$, with coefficient $q - q^{-1}$. For fixed $(k,\ell)$ the only $(a,b)$ with
$R^{ab}_{k\ell} \neq 0$ are $(a,b) = (k,\ell)$, with coefficient $q$ if $k = \ell$ and $1$
otherwise, and $(a,b) = (\ell,k)$ when $k > \ell$, with coefficient $q - q^{-1}$. Hence
\[
  \mathrm{LHS} =
  \begin{cases}
    q\,x_{ik}x_{i\ell}, & i = j,\\[2pt]
    x_{ik}x_{j\ell} + (q-q^{-1})\,x_{jk}x_{i\ell}, & i < j,\\[2pt]
    x_{ik}x_{j\ell}, & i > j,
  \end{cases}
  \qquad
  \mathrm{RHS} =
  \begin{cases}
    q\,x_{j\ell}x_{ik}, & k = \ell,\\[2pt]
    x_{j\ell}x_{ik}, & k < \ell,\\[2pt]
    x_{j\ell}x_{ik} + (q-q^{-1})\,x_{jk}x_{i\ell}, & k > \ell.
  \end{cases}
\]
We record the six essential cases; the remaining ones are these re-indexed.

\emph{(a) $i = j$, $k = \ell$:} both sides are $q\,x_{ik}x_{ik}$; no relation.

\emph{(b) $i = j$, $k < \ell$:} $q\,x_{ik}x_{i\ell} = x_{i\ell}x_{ik}$, which is (R1).

\emph{(c) $i = j$, $k > \ell$:}
$q\,x_{ik}x_{i\ell} = x_{i\ell}x_{ik} + (q-q^{-1})x_{ik}x_{i\ell}$, i.e.
$q^{-1}x_{ik}x_{i\ell} = x_{i\ell}x_{ik}$, which is (R1) with $\ell < k$.

\emph{(d) $i < j$, $k = \ell$:}
$x_{ik}x_{jk} + (q-q^{-1})x_{jk}x_{ik} = q\,x_{jk}x_{ik}$, i.e.
$x_{ik}x_{jk} = q^{-1}x_{jk}x_{ik}$, which is (R2).

\emph{(e) $i < j$, $k > \ell$:} the $(q-q^{-1})$ terms coincide and cancel, leaving
$x_{ik}x_{j\ell} = x_{j\ell}x_{ik}$ with $i<j$ and $\ell<k$, which is (R3).

\emph{(f) $i < j$, $k < \ell$:}
$x_{ik}x_{j\ell} + (q-q^{-1})x_{jk}x_{i\ell} = x_{j\ell}x_{ik}$, i.e.
$x_{j\ell}x_{ik} - x_{ik}x_{j\ell} = (q-q^{-1})x_{jk}x_{i\ell}$; since $x_{jk}$ and $x_{i\ell}$
commute by (e), this is (R4).

Conversely (R1)--(R4) clearly imply all instances of \eqref{eq:RTTentry}.
\end{proof}

\begin{definition}\label{def:Oq}
The quantum matrix algebra $\mathcal{O}_q(M_n)$ is the $\kk$-algebra on generators
$\{t_{ab}\}_{1\le a,b\le n}$ subject to, for all $i<k$ and $j<\ell$,
\[
  t_{i\ell}t_{ij} = q\,t_{ij}t_{i\ell}, \quad
  t_{kj}t_{ij} = q\,t_{ij}t_{kj}, \quad
  t_{i\ell}t_{kj} = t_{kj}t_{i\ell}, \quad
  t_{k\ell}t_{ij} - t_{ij}t_{k\ell} = (q-q^{-1})\,t_{i\ell}t_{kj},
\]
a bialgebra via $\Delta(t_{ab}) = \sum_c t_{ac}\otimes t_{cb}$, $\varepsilon(t_{ab}) = \delta_{ab}$.
\end{definition}

\begin{remark}
The relations of Definition 4.8 are the FRT relations of $\mathcal{O}_q(M_n)$ in
the orientation actually produced by the matrix (3). They are obtained from the
more commonly written presentation
\[
  t_{ij}t_{i\ell} = q\,t_{i\ell}t_{ij},\quad
  t_{ij}t_{kj} = q\,t_{kj}t_{ij},\quad
  t_{ij}t_{k\ell} = t_{k\ell}t_{ij}\ (i<k,\ j>\ell),
\]
\[
  t_{ij}t_{k\ell} - t_{k\ell}t_{ij} = (q - q^{-1})\,t_{i\ell}t_{kj}
  \quad (i<k,\ j<\ell)
\]
by substituting $q \mapsto q^{-1}$. The two presentations define isomorphic
bialgebras. Writing $\sigma(a) := n+1-a$, the assignment
$t_{ab} \mapsto t_{\sigma(a)\sigma(b)}$ carries one set of relations to the other,
since $\sigma$ reverses the order on $\{1,\dots,n\}$, and it is compatible with
$\Delta$ and $\varepsilon$ because $\sum_c t_{ac} \otimes t_{cb}$ is preserved
under the re-indexing $c \mapsto \sigma(c)$. 
\end{remark}

\begin{theorem} \label{O_q}
Let $\mathcal{R}_n$ be the rose quiver with one vertex and $n$ loops, and let $R$ be the matrix given as 
\[R 
  \;=\;
  q \sum_{i=1}^n E_{ii} \otimes E_{ii}
  \;+\;
  \sum_{\substack{i,j=1 \\ i\neq j}}^n E_{ii} \otimes E_{jj}
  \;+\;
  (q - q^{-1}) \sum_{\substack{i,j=1 \\ i<j}}^n E_{ij} \otimes E_{ji}.
\]
Then there is an isomorphism of
bialgebras
\[
 \mathfrak H^{\mathrm{def}}(\mathcal{R}_n) \;\cong\; \mathcal{O}_q(M_n).
\]
\end{theorem}

\begin{proof}
Define $\varphi: \mathfrak H^{\mathrm{def}}(\mathcal R_n)\rightarrow \mathcal O_q(M_n)$ as $\varphi(x_{ab})= t_{ab}$.
By Proposition~\ref{prop:relations} the defining relations of
$\mathfrak H^{\mathrm{def}}(\mathcal{R}_n)$ in the generators $x_{ab}$ are exactly (R1)--(R4), i.e.\ the defining
relations of $\mathcal{O}_q(M_n)$ under $\varphi(x_{ab}) = t_{ab}$; so $\varphi$ is an algebra
isomorphism. The rose $R_n$ is the multi-rose of type $(1,n)$, and $I = I_{\mathrm{loop}}$ in
this case is a biideal (to be seen, more generally in Lemma 4.13 and Lemma 4.14) and consequently, $\mathfrak H^{\mathrm{def}}(R_n)$ is a weak bialgebra. Since $|Q_0| = 1$ the
counital subalgebras are $k$, so it is an ordinary bialgebra. By Proposition 4.7
the defining relations of $\mathfrak H^{\mathrm{def}}(R_n)$ in the generators $x_{ab}$ are
exactly (R1)--(R4), i.e.\ the defining relations of $\mathcal{O}_q(M_n)$; hence
$\varphi(x_{ab}) = t_{ab}$ is an algebra isomorphism. It is a bialgebra map
because the face coproduct of Definition 4.1 specializes at the single vertex
$[v,v]$ to $\Delta(x_{ab}) = \sum_c x_{ac} \otimes x_{cb}$ and
$\varepsilon(x_{ab}) = \delta_{ab}$, which are the structure maps of
Definition 4.8.
\end{proof}

\begin{remark} 
$\mathfrak H(\mathcal{R}_n) = \kk\langle x_{ab} : 1 \le a,b \le n\rangle$ is free on $n^2$ generators, so
$\mathrm{GKdim}\,\mathfrak H(\mathcal{R}_n) = \infty$, whereas
$\mathrm{GKdim}\,\mathcal{O}_q(M_n) = n^2$. Thus $\mathfrak H^{\mathrm{def}}(\mathcal{R}_n)$ is not a flat deformation of
$\mathfrak H(\mathcal{R}_n)$.
\end{remark}

Now, we proceed to generalize this construction. Given $R \in \End(\kk^{Q_1}\otimes \kk^{Q_1})$ with entries $R^{pq}_{rs}$ defined by
$R(e_r \otimes e_s) = \sum_{p,q} R^{pq}_{rs}\, e_p \otimes e_q$, and the full arrow matrix
$X = (x_{p,q})_{p,q \in Q_1}$, one has
\[
  (X_1X_2)^{pq}_{rs} = x_{p,r}\,x_{q,s},
  \qquad
  (X_2X_1)^{pq}_{rs} = x_{q,s}\,x_{p,r},
\]
and, writing $\Phi := RX_1X_2 - X_2X_1R$,
\begin{equation}\label{eq:Phi}
  \Phi^{pq}_{rs} \;=\; \sum_{a,b \in Q_1} R^{pq}_{ab}\; x_{a,r}x_{b,s}
    \;-\; \sum_{a,b \in Q_1} R^{ab}_{rs}\; x_{q,b}x_{p,a}.
\end{equation}

We say $R$ \emph{preserves the source--range bigrading} if it preserves the decomposition
$\kk^{Q_1} \otimes \kk^{Q_1} = \bigoplus_{u,v,w,z} \kk^{Q_1(u,v)} \otimes \kk^{Q_1(w,z)}$,
where $Q_1(u,v) = \{e \in Q_1 : s(e) = u,\ r(e) = v\}$.

\begin{definition}\label{def:multirose}
A quiver $Q$ is called a \emph{multi-rose of type $(m,n)$} if $|Q_0| = m$, every arrow of $Q$ is a loop, and
each vertex $u \in Q_0$ carries exactly $n$ loops $\ell^u_1, \dots, \ell^u_n$. Thus $Q$ is the
disjoint union of $m$ copies of the rose $\mathcal{R}_n$.
\end{definition}

Fix a multi-rose $Q$ of type $(m,n)$ and a Hecke $R$-matrix $R \in \End(\kk^n \otimes \kk^n)$
(for instance the standard $GL_q(n)$ matrix). For $u, w \in Q_0$ let
\[
  X^{[u,w]} \;=\; \big(x_{\ell^u_a,\, \ell^w_b}\big)_{1 \le a,b \le n}
\]
be the corresponding block of the arrow matrix, and set
\begin{equation}\label{eq:Phiuw}
  \Phi[u,w] \;:=\; R\,X^{[u,w]}_1X^{[u,w]}_2 \;-\; X^{[u,w]}_2X^{[u,w]}_1\,R ,
  \qquad
  \Phi[u,w]^{pq}_{cd} = \sum_{a,b} R^{pq}_{ab}\,x_{a,c}x_{b,d} - \sum_{a,b} R^{ab}_{cd}\,x_{q,b}x_{p,a},
\end{equation}
all indices running over $1,\dots,n$ and identified with loops at $u$ (row side) resp.\ $w$
(column side). Let $I$ be the two-sided ideal of $\mathfrak H(Q)$ generated by all entries of all
$\Phi[u,w]$, $u,w \in Q_0$.

Set $H_q(Q) := H(Q)/I$. When $|Q_0| = 1$ every arrow is a loop at the unique
vertex and $I = I_{\mathrm{loop}}$, so $\mathfrak H_q(Q) = \mathfrak H^{\mathrm{def}}(Q)$; for
$|Q_0| \ge 2$ the inclusion $I_{\mathrm{loop}} \subsetneq I$ is strict, and it is
$I$, not $I_{\mathrm{loop}}$, that is a biideal.

\begin{remark}\label{rem:global}
The family $\{\Phi[u,w]\}$ \emph{is} of the global form \eqref{eq:Phi}, for the bigrading-preserving
matrix
\[
  \widetilde R \;:=\; \Big(\bigoplus_{u \in Q_0} R\big|_{\kk^{L_u} \otimes \kk^{L_u}}\Big)
    \;\oplus\; \id\big|_{\bigoplus_{u \neq w} \kk^{L_u}\otimes\kk^{L_w}},
    \qquad L_u := \{\ell^u_1,\dots,\ell^u_n\}.
\]
Indeed, if $p \in L_u$ and $q \in L_{u'}$ with $u \neq u'$, then in \eqref{eq:Phi} the first term
requires $r(p) = s(q)$ and the second requires $r(q) = s(p)$, both impossible, so
$\Phi^{pq}_{rs} = 0$; and if $p,q \in L_u$, $r,s \in L_w$, then $\Phi^{pq}_{rs} = \Phi[u,w]^{pq}_{rs}$.
\end{remark}

\begin{lemma} \label{lem:coideal}
With notation as above, for all $p,q,r,s \in Q_1$,
\[
  \Delta\big(\Phi^{pq}_{rs}\big)
  \;=\; \sum_{e,f \in Q_1} \Phi^{pq}_{ef} \otimes (X_1X_2)^{ef}_{rs}
     \;+\; \sum_{e,f \in Q_1} (X_2X_1)^{pq}_{ef} \otimes \Phi^{ef}_{rs}.
\]
In particular $\Delta(I) \subseteq I \otimes \mathfrak H(Q) + \mathfrak H(Q) \otimes I$.
\end{lemma}

\begin{proof}
$\Delta$ is an algebra map and $\Delta(x_{a,c}) = \sum_{e \in Q_1} x_{a,e} \otimes x_{e,c}$, so
\[
  \Delta(x_{a,c}x_{b,d}) = \sum_{e,f} x_{a,e}x_{b,f} \otimes x_{e,c}x_{f,d},
  \qquad
  \Delta(x_{q,b}x_{p,a}) = \sum_{e,f} x_{q,f}x_{p,e}\otimes x_{f,b}x_{e,a}.
\]
Substituting into \eqref{eq:Phi} and recognising the inner sums,
\[
  \Delta\big(\Phi^{pq}_{rs}\big)
  = \sum_{e,f}\Big[(RX_1X_2)^{pq}_{ef} \otimes (X_1X_2)^{ef}_{rs}
    - (X_2X_1)^{pq}_{ef}\otimes (X_2X_1R)^{ef}_{rs}\Big].
\]
Now add and subtract $\sum_{e,f}(X_2X_1R)^{pq}_{ef}\otimes (X_1X_2)^{ef}_{rs}$; the first
difference is $\sum_{e,f}\Phi^{pq}_{ef}\otimes(X_1X_2)^{ef}_{rs}$, while
\[
  \sum_{e,f}(X_2X_1R)^{pq}_{ef}\otimes(X_1X_2)^{ef}_{rs}
  = \sum_{g,h}(X_2X_1)^{pq}_{gh}\otimes\Big(\sum_{e,f}R^{gh}_{ef}(X_1X_2)^{ef}_{rs}\Big)
  = \sum_{g,h}(X_2X_1)^{pq}_{gh}\otimes (RX_1X_2)^{gh}_{rs},
\]
so the remaining two terms combine into $\sum_{g,h}(X_2X_1)^{pq}_{gh}\otimes\Phi^{gh}_{rs}$.
The last assertion follows because $I \otimes \mathfrak H(Q) + \mathfrak H(Q) \otimes I$ is an ideal of $\mathfrak H(Q) \otimes \mathfrak H(Q)$ and
$\Delta$ is multiplicative.
\end{proof}

Note that the above lemma holds for the face coproduct on the
full arrow matrix of any finite quiver. 

\begin{lemma} \label{lem:counit}
Let $Q$ be a multi-rose of type $(m,n)$ and $\Phi[u,w]$ as in \eqref{eq:Phiuw}. Then
$\varepsilon(\alpha\,\Phi[u,w]^{pq}_{cd}\,\beta) = 0$ for all $\alpha,\beta \in \mathfrak H(Q)$; that is,
$I \subseteq \ker\varepsilon$.
\end{lemma}

\begin{proof}
By linearity it suffices to take $\alpha = x_{g,h}$, $\beta = x_{g',h'}$ with $g,h,g',h'$ paths of
equal length; since all arrows are loops, every path is a power of loops at a single vertex, and
$\varepsilon$ of a basis element $x_{A,B}$ is $\delta_{A,B}$.

Consider the first sum of \eqref{eq:Phiuw}. The product
$x_{g,h}\,x_{a,c}x_{b,d}\,x_{g',h'} = x_{g a b g',\; h c d h'}$ is nonzero only if
$g, a, b, g'$ are all loops at one vertex, necessarily $u$, since $a,b \in L_u$ and
$h, c, d, h'$ are loops at one vertex, necessarily $w$. Its counit is
$\delta_{gabg',\,hcdh'}$, which is $1$ exactly when $g = h$ (forcing $u = w$), $a = c$, $b = d$ and
$g' = h'$. The total contribution of the first sum is therefore
$\delta_{g,h}\delta_{g',h'}\delta_{u,w}\,R^{pq}_{cd}$.

For the second sum, $x_{g,h}\,x_{q,b}x_{p,a}\,x_{g',h'} = x_{gqpg',\,hbah'}$ has counit $1$ exactly
when $g = h$, $q = b$, $p = a$ and $g' = h'$, again forcing $u = w$; the coefficient is
$R^{ab}_{cd} = R^{pq}_{cd}$. So the second sum contributes
$\delta_{g,h}\delta_{g',h'}\delta_{u,w}\,R^{pq}_{cd}$ as well, and the two cancel.
\end{proof}

\begin{theorem}\label{thm:main}
Let $Q$ be a multi-rose of type $(m,n)$ with $m = |Q_0| \ge 1$, $n \ge 1$, let $R$ be a Hecke
$R$-matrix on $\kk^n \otimes \kk^n$, and let $I \subseteq \mathfrak H(Q)$ be the two-sided ideal generated by
the entries of $\Phi[u,w]$ for all $u,w \in Q_0$. Then:
\begin{enumerate}
  \item $I$ is a biideal, so
        $\mathfrak H_q(Q)= \mathfrak H(Q)/I$ is an $\mathbb{N}$-graded weak bialgebra.
  \item Taking $R$ to be the standard $GL_q(n)$ $R$-matrix, there is an isomorphism of algebras
        \[
          \mathfrak H_q(Q) \;\cong\; \bigoplus_{(u,w) \in Q_0 \times Q_0} \mathcal{O}_q(M_n),
        \]
        with comultiplication
        $\Delta\big(x^{[u,w]}_{ab}\big) = \sum_{w' \in Q_0}\sum_{c=1}^{n} x^{[u,w']}_{ac}\otimes x^{[w',w]}_{cb}$
        and counit $\varepsilon\big(x^{[u,w]}_{ab}\big) = \delta_{u,w}\delta_{a,b}$. Equivalently,
        $\mathfrak H^{\mathrm{def}}(Q) \cong \mathcal{C}_{Q_0} \otimes \mathcal{O}_q(M_n)$, where $\mathcal{C}_{Q_0}$ is the
        $m \times m$ comatrix face algebra of the quiver with $m$ vertices and no arrows.
  \item The counital subalgebras satisfy $H_s \cong H_t \cong \kk^{Q_0}$, of dimension $m$. Hence
        Hence $\mathfrak H_q(Q)$ is a bialgebra is a bialgebra if and only if $m = 1$, in which case it is $\mathcal{O}_q(M_n)$
        (recovering Theorem 4.9); for $m \ge 2$ it is a genuine quantum groupoid whose base is
        $m$-dimensional and whose deformation is nontrivial in every one of the $m^2$ corners.
\end{enumerate}
\end{theorem}

\begin{proof}
(1) $I$ is an ideal by construction, a coideal by Lemma~\ref{lem:coideal} together with
Remark~\ref{rem:global}, and lies in $\ker\varepsilon$ by Lemma~\ref{lem:counit}. Since the
generators of $I$ are homogeneous of degree $2$, the grading descends.

(2) By Proposition 4.2, $\mathfrak H(Q) \cong \kk\Qhat$, and since every arrow of $Q$ is a loop, $\Qhat$ has
vertex set $Q_0 \times Q_0$ and exactly $n^2$ loops $x^{[u,w]}_{ab} := x_{\ell^u_a, \ell^w_b}$ at
$[u,w]$, with no other arrows. Hence $\kk\Qhat \cong \bigoplus_{u,w}\kk\langle x^{[u,w]}_{ab}\rangle$
as an algebra, the summands being the pairwise orthogonal corners
$x_{u,w}\,\mathfrak H(Q)\,x_{u,w}$. The relations \eqref{eq:Phiuw} are supported in the corner $[u,w]$ and are,
by Proposition 4.7, exactly the FRT relations of
$\mathcal{O}_q(M_n)$ in the generators $x^{[u,w]}_{ab}$. The formulas for $\Delta$ and
$\varepsilon$ are those of Definition 4.1 read off in this basis.

(3) The counital subalgebras of $\mathfrak H(Q)$ are spanned by the length-zero generators $x_{i,j}$, and
these are untouched by $I$ (which is generated in degree $2$), so they survive in the quotient;
they are isomorphic to $\kk^{Q_0}$. The remaining statements follow.
\end{proof}

A weak bialgebra $H$ is a \emph{universal quantum linear semigroupoid} (UQSGd) of
a graded algebra $A$ if $A$ is an $H$-comodule algebra (with the base of $H$
matching $A_0$) and $H$ is universal with this property. In \cite{HWWW}, left, right,
and transposed versions are defined. 

\begin{theorem}[{\cite[Thm.~1.8]{HWWW}}]\label{thm:hwww18}
For a finite quiver $Q$, the left, right, and transposed UQSGds of the path
algebra $\Bbbk Q$ exist and are each isomorphic, as weak bialgebras, to Hayashi's
face algebra $\mathfrak H(Q)$.
\end{theorem}

This is the precise sense in which the face algebra is the universal coacting object: $\mathfrak H(Q)\cong\Bbbk\widehat Q$ 
\cite{CW} is the UQSGd of $\Bbbk Q$. 

\begin{theorem}[{\cite[Prop.~1.9]{HWWW}}]\label{thm:hwww19}
Let $I\subseteq\Bbbk Q$ be a graded ideal generated in degrees $\ge2$. If the
UQSGd $O_{*}(\Bbbk Q/I)$ exists \textup{(}$*\in\{\mathrm{left},\mathrm{right},\mathrm{trans}\}$\textup{)},
then there is a biideal $\mathcal I\subseteq\mathfrak H(Q)$ and an isomorphism of
weak bialgebras
\[
O_{*}(\Bbbk Q/I)\;\cong\;\mathfrak H(Q)/\mathcal I .
\]
\end{theorem}

\begin{conjecture} \cite[Conj.~1.6]{HWWW} \label{conj:hwww}
Over an algebraically closed field, every finite-dimensional weak bialgebra with
commutative counital subalgebras is isomorphic to a weak bialgebra quotient of
$\mathfrak H(Q)$ for some finite quiver $Q$.
\end{conjecture}

It is the weak-bialgebra analogue of the classical fact that every
finite-dimensional algebra over an algebraically closed field is a quotient of
some path algebra. 

For the rose quiver $\mathcal R_n$, $\Bbbk \mathcal R_n=\Bbbk\langle t_1,\dots,t_n\rangle$
is connected, so its UQSGds reduce to Manin's universal bialgebras \cite{Manin-quantum-groups},
and Theorem~\ref{thm:hwww18} gives
$\mathfrak H(\mathcal R_n)\cong O_{*}(\Bbbk\langle t_1,\dots,t_n\rangle)$.

\begin{remark}
Let $A^n_q := \Bbbk\langle t_1,\dots,t_n\rangle/(t_it_j - q\,t_jt_i : i < j)$ denote
quantum affine space, a degree-two graded quotient of $\Bbbk \mathcal R_n$. For a multi-rose quiver 
$Q$ of type $(m,n)$ one has $\Bbbk Q = \bigoplus_{u \in Q_0} \Bbbk\langle t^u_1,\dots,t^u_n\rangle$,
and the graded quotient $A := \bigoplus_{u \in Q_0} A^n_q$ by the quantum-affine-space
relations in each summand is generated in degree two. Theorem 4.17 therefore
applies and identifies the corresponding universal quantum semigroupoid with a
biideal quotient of $\mathfrak H(Q)$; Theorem 4.15 exhibits that quotient explicitly. For
$m = 1$ this recovers the statement that $\mathcal{O}_q(M_n)$ is the transposed
universal bialgebra of $A^n_q$ in the sense of Manin [9].
\end{remark}

\smallskip

\noindent {\bf Author Contributions.}  CG and AKS together wrote this article.

\smallskip

\noindent {\bf Data Availability.} No datasets were generated or analyzed during the current study.
\smallskip

\noindent {\bf Declarations.}
\\
{\bf Competing interests.} The authors declare no competing interests.\\
{\bf Ethical Approval.} Not applicable.

\end{document}